\documentclass{article}                     
\usepackage{graphicx, natbib}

\usepackage{mathptmx}      
\usepackage{amssymb, amsmath, verbatim, amsfonts, amsthm}

\newtheorem{theorem}{Theorem}
\newtheorem{lemma}{Lemma}

\bibpunct{[}{]}{;}{n}{}{,} 
\renewcommand{\theequation}{\arabic{section}.\arabic{equation}}

\renewcommand{\v}{v}

\newcommand{\vs}{v^*}
\newcommand{\p}{p}

\newcommand{\bv}{\bar{v}}
\newcommand{\bvs}{\bar{v}^*}
\newcommand{\bp}{\bar{p}}
\newcommand{\bdelta}{\bar{\delta}}

\newcommand{\ddbp}{\dot{\bar{\p}}}
\newcommand{\dbv}{\dot{\bar{\v}}}
\newcommand{\dbvs}{\dot{\bar{\v}}^*}

\newcommand{\tv}{\tilde{v}}
\newcommand{\tvs}{\tilde{v}^*}
\newcommand{\tp}{\tilde{p}}

\renewcommand{\u}{u}
\newcommand{\bu}{\bar{u}}

\newcommand{\w}{w}

\newcommand{\V}{\mathbb{V}}
\renewcommand{\P}{\mathbb{P}}
\newcommand{\T}{\mathbb{T}}

\newcommand{\TO}{T_\text{s}}
\newcommand{\Tvd}{T_I}
\newcommand{\Tvb}{T_{II}}
\newcommand{\Tvsd}{T_{III}}
\newcommand{\Tvsb}{T_{IV}}
\newcommand{\TgenI}{T_{\text{genetic},1}}
\newcommand{\TgenII}{T_{\text{genetic},2}}
\newcommand{\Tgen}{T_{\text{genetic}}}
\newcommand{\Twild}{T_\text{wild}}
\newcommand{\Tmut}{T_\text{mutant}}

\newcommand{\bTvd}{\bar{T}_I}
\newcommand{\bTvb}{\bar{T}_{II}}
\newcommand{\bTvsd}{\bar{T}_{III}}
\newcommand{\bTvsb}{\bar{T}_{IV}}

\newcommand{\hht}{\hat{t}}
\newcommand{\se}{s_{II}}
\newcommand{\ses}{s_{IV}}
\renewcommand{\H}{H}
\newcommand{\HII}{\H_{II}}
\newcommand{\HIV}{\H_{IV}}

\newcommand{\kval}{m}

\newcommand{\philim}{\phi_\text{lim}}
\newcommand{\psilim}{\psi_\text{lim}}
\newcommand{\tfinal}{t_f}

\renewcommand{\w}{w}
\newcommand{\z}{\Xi}
\renewcommand{\w}{E}

\begin{document}

\title{Analysis of a Stochastic Predator-Prey Model with Applications to Intrahost HIV Genetic Diversity 
}


\author{Sivan Leviyang\thanks{Georgetown University.  Department of Mathematics.  sr286@georgetown.edu}}


\maketitle

\begin{abstract}
During an infection, HIV experiences strong selection by immune system T cells.  Recent experimental work has shown that MHC escape mutations form an important pathway for HIV to avoid such selection.  In this paper, we study a model of MHC escape mutation.  The model is a predator-prey model with two prey, composed of two HIV variants, and one predator, the immune system CD8 cells.  We assume that one HIV variant is visible to CD8 cells and one is not.  The model takes the form of a system of stochastic differential equations. Motivated by well-known results concerning the short life-cycle of HIV intrahost, we assume that HIV population dynamics occur on a faster time scale then CD8 population dynamics.  This separation of time scales allows us to analyze our model using an asymptotic approach.

Using this model we study the impact of an MHC escape mutation on the population dynamics and genetic evolution of the intrahost HIV population.  From the perspective of population dynamics, we show that the competition between the visible and invisible HIV variants can reach steady states in which either a single variant exists or in which coexistence occurs depending on the parameter regime.  We show that in some parameter regimes the end state of the system is stochastic. From a genetics perspective, we study the impact of the population dynamics on the lineages of HIV samples taken after an escape mutation occurs.  We show that the lineages go through severe bottlenecks and that the lineage distribution can be characterized by a Kingman coalescent.  
\end{abstract}


\section{Introduction} \label{S:introduction}

	During HIV infection, HIV and the immune system T cell populations form a complex and dynamic coupled system. Many authors have modeled and analyzed this interaction, e.g. \cite{Nowak_and_May_Book, Perelson_Nature_Reviews_2002}.  In such work, authors usually consider the interaction of  HIV, CD4 T cells, and CD8 T cells through deterministic predator-prey ODE models.  In this paper, our biological motivation is to consider the effect of MHC escape mutations (described biologically below) on the genetic diversity  and population dynamics of the infecting HIV population.  Mathematically, in order to examine these biological issues,  we extend the typical ODE models by including stochastic effects in the population dynamics and considering lineages of infected HIV cells as one looks backward in time.  More specifically, we change the typical ODE models to stochastic differential equations (SDE) and consider a coalescent process of HIV evolution on top of the population dynamics.     
	
	We consider a model of HIV-CD8 interaction that focuses on so-called MHC escape mutations.  Roughly speaking, when HIV enters a CD4 cell certain mechanisms within the cell cut up HIV proteins into small pieces (usually 8-11 amino acids long) and present these pieces on the surface of the cell.  This presentation is accomplished by the binding of the viral pieces to so-called MHC I molecules to form a peptide-MHC complex (pMHC) \cite{deFranco_book}. For our purposes and to simplify the explanation, the pMHC complex can be thought of as representing a certain short nucleotide sequence in the HIV genome.  CD8 cells are equipped with T cell receptors (TCRs) that can bind to a pMHC complex and then destroy the presenting cell.  Critically, each TCR binds to a limited pattern of nucleotide sequences  \cite{deFranco_book}.  In this sense, since each CD8 cell has only one type of TCR, we can think of each CD8 as targeting some short segment of the viral genome.  
	
	Recently, there has been much experimental and statistical work on mutations in HIV that avoid MHC I presentation (e.g. \cite{Althaus_PLOS_Comp_Bio_2008, Delport_PLOS_Pathogens_2008, Frahm_Nature_Immun_2006, Ngumbela_AIDS_Res_Hum_Retro_2008, Rousseau_J_Virology_2008}) and many authors have suggested that such mutations play a key role in HIV dynamics \cite{Goulder_Nature_Reviews_2008, Leslie_JEM_2005, Bhattacharya_Science_2007}.  The MHC I molecule cannot present every type of nucleotide sequence \cite{deFranco_book}, and it is possible for the virus to mutate and evade MHC presentation.  In such a case, a mutation (or series of mutations) will render the virus invisible to the CD8 cell that was previously able to attack virus infected cells.
	
	We consider a simple model of such an MHC escape mutation.  We assume that initially all HIV infected cells are subject to attack by a collection of CD8 cells with a shared TCR.   We refer to CD4 cells infected by HIV variants that are visible to these CD8 cells as wild type cells.  Then, we assume that a single infected cell changes in its HIV genetic state and becomes invisible to the CD8 cell attack, we refer to this new type of infected cell as the mutant type.  Our motivation centers on understanding the population dynamics and genetic evolution caused by escape mutations that occur during the chronic stage of HIV when CD8 and HIV population sizes are relatively stable \cite{Emini_HIV_Book}.  In terms of population dynamics, we are interested in whether the mutant type survives, the wild type survives, or both.  From a genetics perspective, the escape mutation reflects a change in the HIV genome at a given location, but we are interested in the effect of the resultant population dynamics on other parts of the genome.   As an analogy, strong selective sweeps, for example, may be caused by mutations at a given point in the genome, but other areas of the genome are affected by the sweep  \cite{Hartl_Book_Principles_Population}.  In some sense, our model represents a complex selective sweep involving three players and an uncertain final outcome.  We would like to understand the impact of this complex selective sweep on intrahost HIV genetic diversity.  In this paper we will ignore recombination.
	
	We specify birth and death rates for the CD8, wild type, and mutant type cells and consider the associated coupled birth-death processes.  However, rather than consider the birth-death processes directly, we analyze an associated SDE system.   A critical feature of our analysis is a separation of time scales.  HIV has been shown to evolve on a time scale of hours to days, while the CD8 cells evolve on a time scale of days to weeks \cite{Ho_Nature_1995}.  Technically, this means that the CD8 cells' birth-death process has much lower rates than the wild and mutant type birth-death processes.  We use this separation of time scale to apply an asymptotic analysis to the associated SDE system.  We analyze the SDE system in a large population limit that takes the number of HIV infected cells to infinity.  Since the HIV virion population size is on the order of $10^9$ \cite{Ho_Nature_1995}, this limit is appropriate.
	
	We show that in a certain range of time scale separation, stochastic effects become very important.  Our results demonstrate that in certain parameter regimes one of three events occurs with probability one: the mutant type may be lost, the wild type may be lost, or coexistence may occur.  Further, we show that in other parameter regimes, the probability of each of these three events is strictly between zero and one.  These results contrast with results for the deterministic analogue to our SDE system in which coexistence always occurs.  The degree of the time scale separation connects to the rate of CD8 cell response and so our results show that the end result of MHC escape mutations depends heavily on immune system speed.  
	
	We consider questions of genetics diversity by considering the lineages of a sample of infected cells taken after the escape mutation occurs.  We show that the relationship between the number of samples taken some time after the mutation and the number of corresponding lineages that exist immediately before the mutation can be described through a Kingman coalescent \cite{Durrett_Book_Probability_Models}.  An open question that has received considerable attention is the effective size of HIV during infection.  Effective size is a way of baselining the genetic evolution of a population by comparing it to the genetic evolution of classical Wright-Fisher populations \cite{Durrett_Book_Probability_Models}.  Our results demonstrate that the effective size of HIV is heavily influenced by bottlenecks that occur during MHC escape mutation.  We show that the magnitude of these bottlenecks depends on the rate of CD8 cell response and that in this sense, the effective size of HIV is affected by the speed of the immune system.

\section{The Model}  \label{S:bd}
\setcounter{equation}{0}
\setcounter{lemma}{0}

	We let $\v(t), \vs(t), \p(t)$ be the number of wild type infected cells, mutant type infected cells, and CD8 cells targeting the wild type infected cells respectively (here $\p$ stands for predator).  We assume the dynamics of these populations are given by a birth-death process with the following rates.     

\begin{center}
  \begin{tabular}{| c | c | c | }
    \hline
    type & birth rate per cell & death rate per cell \\ \hline
    $\v$ & $\frac{k}{2} + \frac{\Delta k}{2}$ & $\frac{k}{2} - \frac{\Delta k}{2} + c(\v + \vs) + a\p$  \\ \hline
    $\vs$ & $\frac{k^*}{2} + \frac{\Delta k^*}{2}$ & $\frac{k^*}{2} - \frac{\Delta k^*}{2} + c(\v + \vs)$ \\ \hline
    $\p$ & $\frac{h}{2} + b\v $ & $\frac{h}{2} + d\p$ \\
    \hline
  \end{tabular}
\end{center}

The parameters $k, k^*, \Delta k, \Delta k^*, h$ represent baseline birth and death rates for each of the cell types when interaction between the cell types can be ignored.  Note that we set the birth and death rates of $\p$ equal because CD8 cells require antigenic stimulation to expand in number, on the other hand we can assume that the HIV birth rate exceeds the HIV death rate for both mutant and wild type cells.  The parameter $a$ measures the rate of CD8 cell killing of wild type cells, $c$ is a logistic growth factor representing competition between infected cells for uninfected cells.  Finally $b$ represents the rate of CD8 expansion in the presence of wild type antigen, while $d$ represents a logistic growth factor corresponding to competition between CD8 cells for antigenic stimulation.

\subsection{The Approximating SDE}

In \cite{Kurtz_Book}, Kurtz described the connections between birth-death processes and an approximating SDE.  For HIV, $\v, \vs$ and $\p$ are all of enormous order.  If we rescale the system to makes these variables $O(1)$ we can arrive at the following SDE system which approximates the birth-death processes of $\v, \vs, \p$.  In (\ref{E:SDE}), $\V$ is the order of the infected cell population size and $\v, \vs$ and $\p$ represent the now rescaled population variables (see the appendix for a precise description of the rescaling).
\begin{gather}  \label{E:SDE}
d\v =  \v(1 - (\v + \vs) - \p) dt
	+ \sqrt{\frac{\v(k + (\v + \vs) + \p)}{\V}} dB_1(t),
\\ \notag
d\vs = \vs(f - (\v + \vs)) dt
	+ \sqrt{\frac{\vs(k^* + (\v + \vs))}{\V}} dB_2(t), 
\\ \notag
d\p = \epsilon p(v - \alpha p) dt,
\end{gather}
where $\epsilon = \frac{b}{c}, \alpha = \frac{d}{a}$.  We set $\u(t) = (\v(t), \vs(t), \p(t))$.  

	A theorem of Kurtz, see \cite{Kurtz_Book}, says that such approximations becomes exact as the scaling factor, our $\V$, goes to infinity.  However, our system is not exactly of Kurtz's form as our parameters will be scaled with $\V$, i.e. $\epsilon$ will be taken $O(\frac{1}{\log \V})$.  We have not pursued this technical issue, rather from this point on we take (\ref{E:SDE}) as our description of the evolution of $\u$ and no longer consider the birth death processes.  Although this is an approximation, we believe that our results will hold for the birth-death processes as well.  Further, we point out that the birth-death process is an approximation of underlying dynamics, so to a certain extent considering (\ref{E:SDE}) is no worse than considering the birth-death process.  Finally, notice that there is no stochastic term in the $\p$ equation.  There should be one, but we have dropped it.  This will have no effect on our results and simplifies the explanation (see section \ref{S:overview} for a more precise justification).

	We assume that the escape mutation arises at time $t=0$, and that previously $\vs = 0$ and the system (\ref{E:SDE}), restricted to $\v, \p$, is in equilibrium.  This assumption of equilibrium corresponds to our interest in the chronic stage of HIV infection.  More precisely, we take
\begin{equation}  \label{E:initial_conditions}
\v(0) = \frac{\alpha}{1 + \alpha}, \quad
\vs(0) = \frac{1}{\V}, \quad
\p(0) = \frac{1}{1 + \alpha}
\end{equation}

In (\ref{E:SDE}) the absolute fitness of $\v$, in the absence of CD8 cell effects, is $1$ while the fitness of $\vs$ is $f$.  Typically, in order to avoid immune system attack, escape mutants are less fit than the original wild type and so we take $f < 1$.  We contrast absolute fitness with CD8 influenced fitness.  If CD8 cell attack is considered, the fitness of $\v$ and $\vs$ are $1 - \p$ and $f$ respectively.  If $1 - \p(0) > f$ then the mutant is initially less fit than the wild type and the dynamics are not interesting, indeed the mutant will simply quickly die out.  We restrict our attention to the interesting case of $1 - \p(0) < f$.  
Using (\ref{E:initial_conditions}) this translates to $f - \alpha(1-f) > 0$ and we will assume this condition throughout the rest of this paper.

	We will consider (\ref{E:SDE}) in the limit $\V \to \infty$ with $\epsilon = O(\frac{1}{\log \V})$ over a time interval $[0, O(\frac{1}{\epsilon^2})]$.  The $\epsilon$ scaling will be shown to be the correct scaling to see $O(1)$ stochastic effects, and the time interval scaling is the length of time the system needs to be guaranteed to return to an equilibrium after the escape mutation arises. 

	Associated with (\ref{E:SDE}) is the deterministic analogue in which the stochastic terms are simply dropped.
\begin{gather}  \label{E:vvsp_system}
\dbv = \bv(1 - \bp - \bv - \bvs) \\ \notag
\dbvs = \bvs(f - \bv - \bvs) \\ \notag
\ddbp = \epsilon \bp(\bv - \alpha \bp)
\end{gather} 
with,
\begin{equation}
\bv(0) = \frac{\alpha}{1 + \alpha}, \quad
\bvs(0) = \frac{1}{\V}, \quad
\bp(0) = \frac{1}{1 + \alpha}
\end{equation}
We set $\bu(t) = (\bv(t), \bvs(t), \bp(t))$.  Throughout this paper, we use a bar to distinguish a variable associated with the deterministic system (\ref{E:vvsp_system}) from the corresponding variable associated with the stochastic system (\ref{E:SDE}).  

\subsection{Decomposition of Oscillations}

\begin{figure} [ht] 
\begin{center}
\includegraphics[width=.95\textwidth]{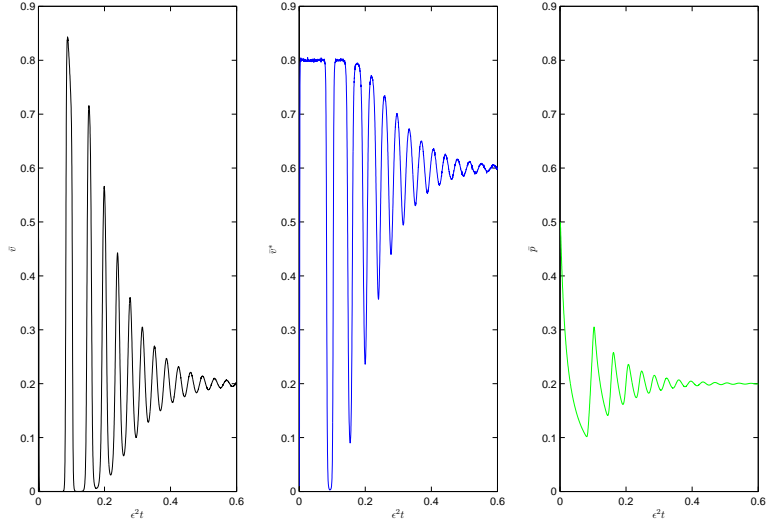}
\caption{Solution of (\ref{E:vvsp_system}) for $\epsilon = .01, \alpha = 1, f = .8$ up to time $O(\frac{1}{\epsilon^2})$}.  
\label{F:full_ODE}
\end{center}
\end{figure}

	The deterministic system (\ref{E:vvsp_system}) produces oscillatory dynamics as the system moves from the original equilibrium of $\bu = (\frac{\alpha}{1 + \alpha}, 0, \frac{1}{1 + \alpha})$ to the new equilibrium of $\bu = (\alpha(1-f),f - \alpha(1-f), 1-f)$.  Figure \ref{F:full_ODE} shows an explicit solution of (\ref{E:vvsp_system}) for $\epsilon = .01, \alpha = 1, f = .8$.  
	
\begin{figure} [ht] 
\begin{center}
\includegraphics[width=.95\textwidth]{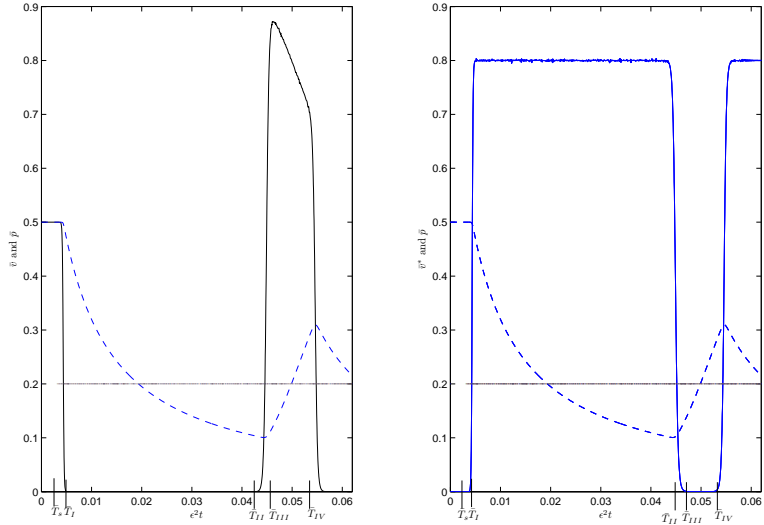}
\caption{Solution of (\ref{E:vvsp_system}) for $\epsilon = .005, \alpha = 1, f = .8$ showing only the first cycle in the dynamics of $\bu$}.  
\label{F:cycle}
\end{center}
\end{figure}

	The oscillatory dynamics result from the time scale separation.  Figure \ref{F:cycle} provides a zoom in of a first oscillation in the dynamics of $\bu$.  Immediately after the escape mutation arises, at $t=0$, there is an initial stage of the dynamics in which $\bv, \bp$ are held relatively fixed while $\bvs$ rises.  This initial stage, which in Figure \ref{F:cycle} ends at $t = \TO$, occurs only immediately after $t=0$ and is not repeated.   After the initial stage ends at time $\TO$, we separate the oscillation into four stages which form a full cycle:  stage I delimited by $[\TO, \bTvd]$, stage II delimited by $[\bTvd, \bTvb]$, stage III delimited by $[\bTvb, \bTvsd]$, and stage IV delimited by $[\bTvsd, \bTvsb]$.  In Stage I, $\bv$ collapses to $o(1)$ levels while $\bvs$ rises to $O(1)$ levels.  This occurs because the CD8 influenced fitness of $\bv$, given by $1 - \bp$, is less than $f$ during Stage I.  In Stage II $\bv$ stays at $o(1)$ levels, but $\bp$ drops until $\bv$ becomes more fit than $\bvs$.  When $\bp = 1-f$, $\bv$ and $\bvs$ are equally fit, indeed at this point $\bv$ reaches its minimum. After that time $\bv$ rises until $\bv(\bTvb) = \bv(\bTvd)$.  In Stage III, $\bvs$ is less fit than $\bv$ and $\bvs$ collapses to $o(1)$ levels.  In Stage IV, $\bp$ rises until $\bvs$ becomes more fit than $\bv$ which causes $\bvs$ to rise until $\bvs(\bTvsb) = \bvs(\bTvsd)$.  At the end of Stage IV the system has returned to the situation of time $\TO$ and the cycle repeats. We refer to the Stages I-IV as a cycle, and so the dynamics of $\u$ and $\bu$ are formed by a sequence of cycles.   As Figure \ref{F:full_ODE} shows, with each cycle the strength of the oscillations is damped.  The stochastic system (\ref{E:SDE}) has identical stages with $\Tvd, \Tvb, \Tvsd, \Tvb$ defined analogously to $\bTvd, \bTvb, \bTvsd, \bTvb$.  $\TO$ is used for both (\ref{E:SDE}) and (\ref{E:vvsp_system}).  The definitions of $\TO, \bTvd, \bTvb, \bTvsd, \bTvb$ will be made precise in section \ref{S:overview}, for now we simply provide the reader with an intuition for the dynamics.

	We will show that (\ref{E:SDE}) behaves essentially as (\ref{E:vvsp_system}).  The exception will be during times in a small subinterval in each of stages II and IV during which $\v$ and $\vs$ are at a minimum.  If $\epsilon$ is too small, then at these subintervals $\v$ or $\vs$ will be driven to zero by the stochastic terms in (\ref{E:SDE}).  If $\epsilon$ is too large, then at these subintervals, (\ref{E:SDE}) will behave as (\ref{E:vvsp_system}) and stochastic effects can be ignored.  However, if $\epsilon$ is at the appropriate scaling, namely $O(\frac{1}{\log \V})$, then we will show that $\v$, in stage II, and $\vs$, in stage IV, behave like Feller diffusions.  
	
	The damping of the oscillations will be crucial to our analysis.  Stochastic effects play an important role when $\v$ and $\vs$ are at their minima.  The damping of the oscillations means that the minima become less extreme with each passing cycle.  Consequently, as we shall show, if $\v$ or $\vs$ are not lost during the first cycle, then (\ref{E:SDE}) will reduce to (\ref{E:vvsp_system}) and stochastic effects can be ignored for the rest of the considered time interval.
	
	From a genetic perspective, the behavior of $\v$ and $\vs$ during Stages II and IV correspond to severe bottlenecks.  We will explore the effects of these bottlenecks on the lineages of a set of samples taken at a time $O(\frac{1}{\epsilon^2})$ time units after the escape mutation occurs.

\section{Results}  \label{S:results}
\setcounter{equation}{0}
\setcounter{lemma}{0}

	We are interested in determining the probability of wild type loss, mutant type loss, or coexistence of the two types.  If the wild type is lost, then the steady state of (\ref{E:SDE}) is given by $\u_M = (0,f,0)$.  If the mutant is lost then the steady state is $\u_W = (\frac{\alpha}{1+\alpha}, 0, \frac{1}{1 + \alpha})$.  And finally if coexistence occurs then the steady state is $\u_C = (\alpha(1-f),f - \alpha(1-f), 1-f)$.

	We further discern between two types of mutant loss.  The mutant may be lost immediately after the initial mutation occurs, that is before the mutant population reaches a significant proportion of the virus population (this will be made precise shortly).   We refer to this as failed mutant dynamics.  Or, the mutant may rise to significant population levels but then subsequently be lost.  We refer to this as the lost mutant dynamics.
	
	Our first result characterizes the probability of these events occurring in the limit of $\V \to \infty$ and the scaling $\epsilon = O(\frac{1}{\log \V})$.  We consider the system up to time $t_f$ such that $t_f = O(\frac{1}{\epsilon^2})$ and we show that by time $t_f$ the system is arbitrarily close to a steady state.  Below, when we write $P(\u \leadsto (u_1, u_2, u_3)) = p$ we mean that for any fixed constant $c$, we have $\lim_{\V \to \infty} P(\|\u(\tfinal) - (u_1,u_2,u_3)\|_\infty < c) = p$.
	
	In order to state our result we need the following definition,
\begin{equation}
\philim(\alpha, f) = \frac{1}{\alpha} \bigg[-(f - \alpha(1-f)) + \log(\frac{1}{(1+\alpha)(1-f))}\bigg].
\end{equation}
\begin{equation}
\psilim(\alpha, f) = \frac{-1}{1 + \alpha} \log(\frac{\alpha H}{(1 + \alpha(H+1))(f - \alpha(1-f))})
		+ (1-f) \log(\frac{(1-f)\alpha \H}{f - \alpha(1-f)}).
\end{equation}
where $\H$ solves the following equality,
\begin{equation}
(1-f)\H = \frac{1}{\alpha} \log(1 + \frac{\alpha}{1 + \alpha} \H).
\end{equation}
Notice that $\H$ is a function of $f$ and $\alpha$ although we do not make this dependence explicit.

\begin{theorem}  \label{T:warmup}
Set $\epsilon = \frac{\beta}{\log \V}$ and  $\tfinal = \frac{t}{\epsilon^2}$ where $\beta$ and $t$ are held fixed.   By a failed mutant we will mean, $\sup_{t' \le t_f} \vs(t') < \epsilon$.  

\vspace{.1cm}
\begin{tabular}{|c||c|}
\hline 
If &  Then in the limit $\V \to \infty$  \\
\hline 
$\frac{\philim(\alpha,f)}{\beta} \ge 1$ 
		& $P(\u \leadsto \u_M) = 1 - p_\text{failed}$
\\ 
& $P(\u \leadsto \u_W) = p_\text{failed}$
\\ \hline 
$\frac{\philim(\alpha,f)}{\beta} < 1$ and $\frac{\psilim(\alpha,f)}{\beta} \ge 1$ 
		& $P(\u \leadsto \u_W, \text{failed mutant}) = p_\text{failed}$ 
\\
		& $P(\u \leadsto \u_W, \text{lost mutant}) = 1 - p_\text{failed}$ 
\\ \hline 
$\frac{\philim(\alpha,f)}{\beta} < 1$ and $\frac{\psilim(\alpha,f)}{\beta} < 1$
		& $P(\u \leadsto \u_C) = 1 - p_\text{failed}$
\\
  & $P(\u \leadsto \u_W) = p_\text{failed}$
\\ \hline
\end{tabular}
\\
where, 
\begin{equation}
p_\text{failed} = \exp[-\frac{4(f - \alpha(1-f))}{(k^*+1)(\alpha+1)}]
\end{equation}
\end{theorem}

\begin{figure} [ht] 
\begin{center}
\includegraphics[width=.95\textwidth]{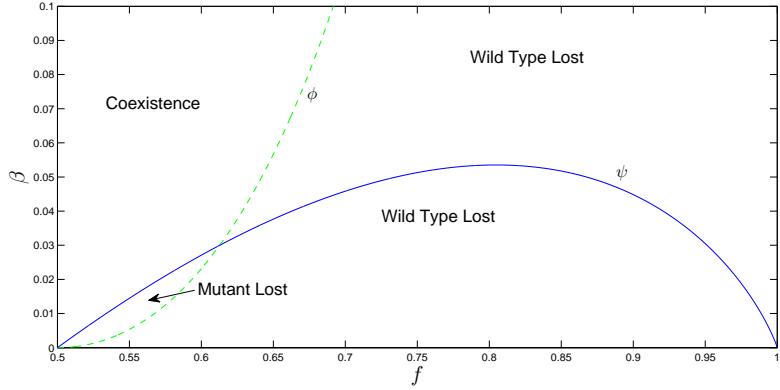}
\caption{Graph of $\philim$ (dashed line) and $\psilim$ (solid line) along with results from the table in Theorem \ref{T:warmup}.  $\philim$ is cut off to allow for presentable scales.  $\alpha = 1$.}.  
\label{F:f_beta}
\end{center}
\end{figure}

	Figure \ref{F:f_beta} graphically displays the results of Theorem \ref{T:warmup} (the case of a failed mutant is ignored).  The relationship between $\philim$ and $\psilim$ seen in Figure \ref{F:f_beta} is general with $\philim$ always increasing and diverging at $f = 1$ and $\psilim$ always possessing a single maximum and equaling zero at the two possible endpoints of $f$:  $\frac{\alpha}{1 + \alpha}$ and $1$.  We see that in the scaling of Theorem \ref{T:warmup}, the end state of (\ref{E:SDE}) is completely deterministic.
	
	If we allow $\beta$ to scale with $\V$ we can analyze the boundaries of the results given in Theorem \ref{T:warmup}.  Define $w(t)$ to be the solution of the following Feller diffusion at time $t$:
\begin{gather}  \label{E:base_Feller}
dw = \sqrt{w} dB. \\ \notag
w(0) = 1.
\end{gather}
Let $\hat{f}(\alpha)$ be the solution to 
\begin{equation}
\philim(\alpha, \hat{f}) = \psilim(\alpha, \hat{f}).
\end{equation}
	
\begin{theorem} \label{T:main}
Set $f = \hat{f}$ and $\philim = \philim(\alpha, f)$.  Set,
\begin{equation}
\epsilon = \philim \bigg(\frac{1}{\log \V} + \frac{\frac{1}{2} \log \log \V}{(\log \V)^2} -
	\frac{\log (\kappa\sqrt{\philim})}{(\log \V)^2}\bigg),
\end{equation}
and $\tfinal = \frac{t}{\epsilon^2}$.  Define,
\begin{gather}
\Twild = \sqrt{\frac{2\pi }{\alpha (1-f)^2}}(\frac{\alpha}{1 + \alpha}) (k+1)\kappa, \\ \notag
\Tmut = \sqrt{\frac{2\pi}{(1-f)(f-\alpha(1-f))}}(k^*+f)(\frac{1}{f})\bigg(\frac{1 + \alpha(1+H)}{(1+\alpha)\alpha(1+H)}\bigg)^H (\eta_{IV})^\frac{\alpha}{H} \kappa.
\end{gather}
where $\eta_{IV}$ is a random variable with distribution,
\begin{equation}  \label{E:eta_IV_def}
\eta_{IV} = w[\Twild].
\end{equation}
Further set,
\begin{gather}
\rho_W = P(w[\Twild] = 0), \\ \notag
\rho_M = P(w[\Tmut] = 0)
\end{gather}
Then in the limit of $\V \to \infty$, 
\begin{gather}
P(\u \leadsto \u_W, \text{failed mutant}) = p_\text{failed} \\ \notag
P(\u \leadsto \u_W, \text{mutant lost}) = (1 - p_\text{failed}) (1 - \rho_W) \rho_M
\\ \notag
P(\u \leadsto \u_M) = (1 - p_\text{failed}) \rho_W
\\ \notag
P(\u \leadsto \u_C) = (1 - p_\text{failed}) (1 - \rho_W) (1 - \rho_M)
\end{gather}
\end{theorem}

	Our techniques allow us to derive similar results for the cases $f < \hat{f}$ and $f > \hat{f}$.  Notice that in Theorem \ref{T:main}, the end state of (\ref{E:SDE}) is truly stochastic.  The scaling of $\epsilon$ in Theorem \ref{T:main} shows that the stochastic regime of (\ref{E:SDE}) scales as $O(\frac{1}{\log \V})$ about $\phi$.  And indeed, the same can be shown to be true about $\psi$.  
	
	We express our genetic results within the context of population sampling.  To make things somewhat concrete, we assume that a sample of $n$ infected cells is taken at time $t_f$.   The genetic composition of this sample of this sample across the HIV genome can be determined if one knows the lineages formed by these $n$ samples and the mutations that occur on these lineages.  The lineages and mutations associated with a sample are not deterministic and hence must be specified through a probability distribution.  
	
	In this work, we do not determine the full distribution of the lineages and mutations associated with the $n$ samples (although our methods should allow for this).  Rather, we characterize the state of the lineages formed by the $n$ samples at time $0$.  To explain this precisely, let $y_1, y_2, \dots, y_n$ be labels for the $n$ infected cells sampled at time $t_f$.  At time $0$ these $n$ samples will have some number of ancestors, say $n_0$, and we can arbitrarily label these ancestors $z_1, z_2, \dots, z_{n_0}$.  Each $z_i$ will be the ancestor of a certain number of the $n$ sampled cells, let this number be $B_i$.  Then we have $B_1 + B_2 + \dots + B_{n_0} = n$ since every sample $y_i$ must be descendant from some $z_j$.  Our results specify the distribution of $(n_0, B_1, B_2, \dots, B_{n_0})$.  
	
	The Kingman coalescent is a specific probability distribution for the lineages of a set of samples corresponding to the lineage distribution in a classical Wright-Fisher population \cite{Durrett_Book_Probability_Models}.  For any time $t$, we let $\Pi(t; n)$ be the lineage distribution specied by the Kingman coalescent $t$ time units prior to sampling of $n$ individuals.  The distribution of $\Pi(t;n)$ is well understood \cite{Durrett_Book_Probability_Models, Tavare_2001_St_Flour_Lecture_Notes}.   Our results will show that the distribution of $(n_0, B_1, B_2, \dots, B_{n_0})$ is the same as that specified by the Kingman coalescent run for some time $\Tgen$.  In this way, we show that our viral population run for time $t_f$ will have a similar level of genetic diversity as a Wright-Fisher population run for time $\Tgen$.   
	
	Our results will actually split into cases depending on whether the wild type is lost, mutant type is lost, or coexistence occurs.  In the case of one type surviving, we will be able to characterize the genetic diversity by a single Kingman coalescent started with $n$ individuals.  However, when the types coexist, we need two Kingman coalescents to characterize the resultant distribution.
		
	The genetic results follow the same scaling regimes as specified in Theorems \ref{T:warmup} and \ref{T:main}.  However, for clarity of presentation we simply state the results associated with the parameter regime given in Theorem \ref{T:main}.  The following theorem shows that the distribution of $(n_0, B_1, B_2, \dots, B_n)$ is a function of four random variables defined in the theorem: $w_1, w_2, \zeta, \xi$. We state the Theorem and then help the reader parse the results.
	
\begin{theorem} \label{T:genetics}
Let $f, \epsilon, \Twild, \Tmut, p_\text{failed}$ be as in Theorem \ref{T:main}.   Let $w_1, w_2$ be independent versions of the Feller process $w$ specified in (\ref{E:base_Feller}) and let the Feller process in the definition of $\Tmut$ be $w_1$. Let $\zeta$ be a uniform random variable on $[0,1]$ and let $\xi$ be a binomial random variable with $n$ trails and success probability $\frac{\alpha (1-f)}{f}$.  Define the times $\TgenI, \TgenII$ according to the following table (in the cases where the entry for $\TgenII$ is a $-$, only $\TgenI$ is defined.)

\begin{tabular}{|c||c|}
\hline 
If &  $\TgenI$, $\TgenII$
\\ \hline 
$\zeta < p_\text{failed}$ & $0$, $-$
\\ \hline
$\zeta > p_\text{failed},  w_1[\Twild] = 0$ & $\infty$, $-$
\\ \hline 
$\zeta > p_\text{failed}, w_1[\Twild] \ne 0, w_2[\Tmut] = 0$ & 
$\Upsilon \int_0^{\sqrt{2\pi}}
							ds \frac{1}{w_1[\Upsilon (k+1) s]}$, $-$
\\ \hline
$\zeta > p_\text{failed}, w_1[\Twild] \ne 0, w_2[\Tmut] \ne 0$
	& $\Upsilon \int_0^{\sqrt{2\pi}}
							ds \frac{1}{w_1[\Upsilon (k+1) s]}$, $\infty$
\\ \hline
\end{tabular}

where $\Upsilon = \sqrt{\frac{1}{\alpha(1-f)^2}}(\frac{\alpha}{1 + \alpha})\kappa$.  In the cases where $\TgenII$ is not defined then,
\begin{equation}
\lim_{\V \to \infty} (n_0, B_1, B_2, \dots, B_{n_0}) = \Pi(\TgenI; n) 
\end{equation}
In the case where $\TgenII$ is defined then,
\begin{equation}
\lim_{\V \to \infty} (n_0, B_1, B_2, \dots, B_{n_0}) = (\Pi(\TgenI; \xi), \Pi(\TgenII, n - \xi)).
\end{equation}
Both of the above limits are meant in the sense of convergence in distribution.
\end{theorem}

	To explain Theorem \ref{T:genetics} we explain the reasoning behind the table.  Take first the case $\zeta < p_\text{failed}$.  In this setting, the mutant fails to enter the viral population and essentially the population stays in $\v, \p$ equilibrium.  No bottlenecks occur and so the $n$ lineages, run for time $O(\frac{1}{\epsilon^2})$ in a population of size $O(\V)$ will not converge.  This is the same as running the Kingman coalescent for zero time.  In the second case, $\zeta > p_\text{failed},  w_1[\Twild] = 0$, the wild type will be lost.  In this case all lineages must come from the original escape mutant.  This is the same as running the Kingman coalescent for infinite time to guarantee convergence of the lineages to a common ancestor.  The third case, $\zeta > p_\text{failed}, w_1[\Twild] \ne 0, w_2[\Twild] = 0$, corresponds to the loss of the mutant type.  A bottle neck occurs in the wild types during Stage II and the mutant is lost during Stage IV.  The bottleneck has an effect on the lineages equivalent to running a Kingman coalescent for $O(1)$ time.  The expression for $\TgenI$ gives the precise length of time the Kingman coalescent needs to be run, notice that $\TgenI$ is stochastic.  Finally, in the case $\zeta > p_\text{failed}, w_1[\Twild] \ne 0, w_2[\Twild] \ne 0$, both the wild and mutant types survive.  In the new equilibrium, the fraction of wild type cells will be $\frac{\alpha (1-f)}{f}$, so if we randomly sample at time $t_f$, $\xi$ gives the probability distribution of the number of wild types sampled.  Working backward in time, all mutant samples come from the original mutant and hence must coalesce by time zero.  $\TgenII$ reflects this, we run the Kingman coalescent infinite time for the part of the sample corresponding to mutant samples.  $\TgenI$ represents wild type samples, and the bottleneck through which they pass in Stage II has an effect on the lineages equivalent to running a Kingman coalescent for time $\TgenI$.  
\section{Discussion} \label{S:discussion}

	In this paper, through a specific model of MHC I escape mutation, we have attempted to emphasize the important interaction between population dynamics and genetic evolution. Indeed, as we have demonstrated, the initial escape mutation brings about oscillatory dynamics in the population sizes of the wild and mutant type cells.  These oscillation then impact the lineages of sampled cells, implying that selection by CD8 cells will affect genetic diversity over the whole viral genome.
	
	Several authors have explored the issue of intrahost HIV effective population size \cite{Brown_PNAS_1997, Rouzine_PNAS_1999, Kuoyos_Trends_MicroBio_2006}. Effective population size, as opposed to census population size, is a measure of the potential genetic diversity of a population.  High effective population sizes imply, at least in the presence of significant mutation, high levels of genetic diversity.   The HIV population is enormous and so one expects the effective population size to be large and, in turn, genetic diversity to be large.  However, experimental data seems to suggest a small level of genetic diversity and hence a small effective population size \cite{Brown_PNAS_1997}.  In this work, we have demonstrated that the notion of an effective population size for HIV in the presence of the predator-prey dynamics of CD8-HIV interaction is subtle (see \cite{Kuoyos_Trends_MicroBio_2006} for a qualitative discussion with similar observations).  The intuition that effective population size and census population size are directly related is an artifact of considering classical Wright-Fisher populations, and HIV is certainly not such a population.  Our results show that significant bottlenecks occur in the HIV population due to CD8 attack, and that the size and time lengths of these bottlenecks is determined not by the population size of HIV but by the rate of CD8 cell response.  In our model $\epsilon$ plays the role of CD8 response rate and $\V$ plays the role of HIV census population size.  We find that if $\epsilon \gg \log \V$ then lineages do not coalesce over time $O(\frac{1}{\epsilon^2})$.  And, in fact, one can show that lineages will coalesce in time $O(\V)$ which leads to an effective size of $O(\V)$.  However if $\epsilon \le O(\log \V)$ then lineages coalesce in time less than $O(\frac{1}{\epsilon^2})$.  In fact, one can show that the bottlenecks last for $O(\frac{1}{\sqrt{\epsilon}})$ time which leads to an effective size of $O(\frac{1}{\sqrt{\epsilon}})$.  
	
	In this paper we have not considered the mutation rate. While effective population size gives some idea of a population's genetic diversity, to fully model genetic evolution one needs to construct a model of lineage distribution, that is a coalescent, and analyze mutations on that lineage distribution.  The strength of mutation then plays a central role in the pattern of genetic diversity.   Data suggests that escape mutations fix or are lost in a matter of weeks.  For a sample of size $n$ the expected number of mutations that take place before the HIV population returns to equilibrium will be bounded by $O(\frac{\mu}{\epsilon^2})$ where $\mu$ is the mutation rate of HIV.  $\mu$ for a single nucleotide is on the order $10^{-5}$, and in this setting we can safely assume that no mutations occur during the escape mutation.  In such a case our results give a complete answer to the effect of our model on genetic diversity since we do not need to know what happens between the time of the escape mutation and the return to equilibrium.  However, if one considers the whole genome, then $\mu$ is on the order of $10^{-1}$ and we must account for mutations.  Our approach should allow for such an analysis, but in this paper we have not pursued this important issue.
	
	The genetic diversity of HIV over the course of the chronic phase of infection has been shown to initially rise, eventually level off, and then finally drop \cite{Shankarappa_J_Virol_1999}.  The drop of HIV genetic diversity is correlated with the onset of AIDS.  We speculate that one possible explanation for this evolution of diversity involves the collapse of the immune system.  Our model suggests that a fast immune system will allow for more diversity.  Indeed, when our $\epsilon$ is large, Theorem \ref{T:genetics} with $\kappa \to 0$ shows that the lineages do not coalesce and so genetic diversity increases with time.    But, as CD4 cells counts collapse, we speculate that the rate of CD8 response will also collapse and $\epsilon$ will drop.  From Theorem \ref{T:genetics} we see, by taking $\kappa \to \infty$, that a small $\epsilon$ forces lineages to coalesce and genetic diversity will drop. This line of reasoning is, of course, highly speculative since the loss of CD4 cells will impact CD8 response and HIV evolution in many ways that our model does not consider.  And of course we are considering CD8 attack at only a single epitope while HIV is typically subject to attack at multiple epitopes.

\section{Proofs of Theorems}  \label{S:overview}
\setcounter{equation}{0}
\setcounter{lemma}{0}

	In this section we prove Theorems \ref{T:warmup}, \ref{T:main}, and \ref{T:genetics}.  As we mentioned in section \ref{S:bd}, we analyze the stochastic system (\ref{E:SDE}) by comparing it to deterministic system (\ref{E:vvsp_system}).  More precisely, in section \ref{S:deterministic} we consider the deterministic system in Stages I-IV.  For each stage, we determine an asymptotic expansion in $\epsilon$ of $\u$ at the end of the stage.  So, for example, in stage II we assume $\bu(\bTvd)$ as given and determine an expansion for $\bu(\bTvb)$.  In section \ref{S:stochastic} we consider the stochastic system through Stage I-IV.  Roughly, the idea will be that in Stage I and III, $\u$ and $\bu$ behave, with high probability, almost identically while in Stage II and IV, $\u$ and $\bu$ differ due to stochastic effects. 
	
	The main technical ideas needed to prove our theorems are contained in the lemmas of section \ref{S:deterministic} and \ref{S:stochastic}.  In this section we put these lemmas together to prove the theorems.  By doing so, we hope to provide the reader with an underlying intuition for the more technical lemmas found in sections \ref{S:deterministic} and \ref{S:stochastic}.  Before beginning in this task, we need to precisely define Stage I-IV and the initial stage.  We set for the deterministic system,
\begin{gather}
\bTvd = \inf\{t>0 : \bv(t) = \epsilon^q\}, \\ \notag
\bTvb = \inf\{t>\bTvd : \bv(t) = \epsilon^q\}, \\ \notag
\bTvsd = \inf\{t>\bTvb : \bvs(t) = \epsilon^q\}, \\ \notag
\bTvsb = \inf\{t>\bTvsd : \bvs(t) = \epsilon^q\},
\end{gather} 
and similarly for the stochastic system, 
\begin{gather}
\Tvd = \inf\{t>0 : \bv(t) = \epsilon^q\}, \\ \notag
\Tvb = \inf\{t>\Tvd : \v(t) = \epsilon^q\}, \\ \notag
\Tvsd = \inf\{t>\Tvb : \vs(t) = \epsilon^q\}, \\ \notag
\Tvsb = \inf\{t>\Tvsd : \vs(t) = \epsilon^q\},
\end{gather}
where $q = 4$.  In section \ref{S:bd} we stated that in section I, $\v$ collapses to $o(1)$ levels.  What we meant, as seen from the definitions above, is that $\v$ collapses to the value of $\epsilon^q$.   Essentially, $q$ is chosen so that when a variable falls below $\epsilon^q$, its effect will not be felt in the $\V \to \infty$ limit.  For example, in Stage II, $\v \le \epsilon^q$ and the result is that the dynamics of $\p$, given by $\dot{\p} = \epsilon \p(\v - \alpha \p)$ can be reduced to the integrable $\dot{\p} = -\epsilon \alpha \p^2$.  Set,
\begin{equation}  \label{E:TO_def}
\TO = \inf\{t : \vs \notin (0,\epsilon^q)\}
\end{equation}
$\TO$ is considered only in the stochastic system and so $\bar{\TO} = \TO$.  Indeed, we will simply start the deterministic system at time $\TO$ by setting $\bu(\TO) = \u(\TO)$.  

	We now proceed to prove Theorems \ref{T:warmup}-\ref{T:genetics}.  Theorems \ref{T:warmup} and \ref{T:main} have similar proofs, so we prove Theorem \ref{T:main} and within that proof we comment on the connections to the proof of Theorem \ref{T:warmup}
	
\begin{proof}[Theorem \ref{T:main}]
To help the reader parse our arguments, we decompose this proof according to which stage we consider.  At the end of the proof we consider a full cycle.  

\vspace{.1cm}
\begin{flushleft}
\textbf{Initial Stage:}\\
\end{flushleft}
We start by considering the initial stage.  From (\ref{E:TO_def}) we have three possibilities.  $\v(\TO) = \epsilon^q$, $\v(\TO) = 0$ or $\TO = \infty$.  Lemma \ref{L:I_0} shows that $P(\TO = \infty) = 0$ and $\lim_{\V \to \infty} P(\v(\TO) = 0) = p_\text{failed}$.  So with probability $p_\text{failed}$ the stochastic system never exceeds $\epsilon^q$ and we have a failed mutant.  Lemma \ref{L:I_0} also shows that $\v, \p$ do not deviate beyond $O(\epsilon^q)$ from $\v(0), \p(0)$ by time $\TO$, so if a failed mutant occurs, the system returns to state $\u_W$.  

\vspace{.1cm}
\begin{flushleft}
\textbf{Stage I:}\\
\end{flushleft}
If the mutant does not fail, then Stage I starts.  We take $\bu(\TO) = \u(\TO)$.    Lemma \ref{L:I_and_III_s} gives for Stage I,
\begin{equation}  \label{E:I_det}
P(\sup_{\TO \le t \le \Tvd} \|u(t) - \bar{u}(t)\|_\infty \ge \frac{1}{\V^\frac{1}{8}})
	\le O(\frac{(\log \log \V)^2 (\log \V)^2}{\V^\frac{1}{4}})
\end{equation}
and Lemma \ref{L:I} gives the asymptotic of $\bu$, up to $O(\epsilon)$, at the end of stage I. 

\vspace{.1cm}
\begin{flushleft}
\textbf{Stage II:}\\
\end{flushleft}
Lemma \ref{L:II} shows that for the deterministic system $\v$ reaches an absolute minimum in Stage II.  In section \ref{S:stochastic} we label that minimum as $\frac{\se}{\epsilon}$ and we are able to asymptotically compute its value as a function of $\p(\Tvd)$.  To emphasize this we write,  $\frac{\se(\p(\Tvd))}{\epsilon}$.  In Lemmas \ref{L:stochastic_II_part1}-\ref{L:stochastic_II_part3} we define an interval $[t_0, t_1]$ centered about  $\frac{\se(\p(\Tvd))}{\epsilon}$ and of width $\frac{2}{\epsilon^\kval}$.  For our proofs to work, $\kval$ can be taken as any value between $\frac{1}{2}$ and $\frac{2}{3}$, but intuitively what matters is that the width of $[t_0, t_1]$ be much greater than $\frac{2}{\sqrt{\epsilon}}$.  Lemma \ref{L:stochastic_II_part1} shows that on $[\Tvd, t_0]$, $|\u - \bu|$ is bounded by $O(\epsilon^\frac{q+1}{2})$ with probability $O(\epsilon^2)$.  Lemma \ref{L:stochastic_II_part3} shows the identical conclusion for the interval $[t_1, \Tvb]$.  Lemma \ref{L:stochastic_II_part2} shows that inside the interval $[t_0, t_1]$ stochastic effects matter and $\u$ may deviate from $\bu$.  Indeed, Lemma \ref{L:stochastic_II_part2} gives,
\begin{equation}  \label{E:P_base_4}
\lim_{\V \to \infty} P(\v \text{ is lost in }[t_0,t_1]) = \lim_{\V \to \infty} P(w\big[\sqrt{2\pi}(k+1)\z_{II}\big] = 0)
\end{equation}
where,
\begin{equation}  \label{E:z_II_warmup_0}
\z_{II} = \sqrt{\frac{1}{\alpha (1-f)^2}} (\frac{\exp[-\phi(\frac{\se(\p(\Tvd))}{\epsilon})]}{\V \epsilon^q \sqrt{\epsilon}})
\end{equation}
The function $\phi$ in (\ref{E:z_II_warmup_0}) is defined in section \ref{S:stochastic}.  We caution the reader that $\philim$ is not $\phi$, essentially $\phi$ is  $\log(\frac{\v(t)}{\v(\Tvd)})$ although we define $\phi$ somewhat differently for technical reasons.  The connection between $\phi$ and $\philim$ is given by the following relation which is justified in Lemma \ref{L:stochastic_II_part2}.
\begin{equation}  \label{E:phi_philim_relation}
-\phi(\frac{\se}{\epsilon}) = \frac{\philim}{\epsilon} - q |\log \epsilon| + \log(\frac{\alpha}{1 + \alpha}) + O(\epsilon^\frac{3}{2} |\log \epsilon|).
\end{equation}
This relation leads to,
\begin{equation}  \label{E:phi_to_philim}
\exp[-\phi(\frac{\se(\p(\Tvd))}{\epsilon})] =  \exp[\frac{\philim}{\epsilon}] \epsilon^q (\frac{\alpha}{1 + \alpha})(1 + O(\epsilon^\frac{3}{2} |\log \epsilon|)).
\end{equation}
Plugging (\ref{E:phi_to_philim}) into (\ref{E:z_II_warmup_0}) and using the scaling of $\epsilon$ in Theorem \ref{T:main}  gives,
\begin{equation}  \label{E:z_II_warmup}
\lim_{\V \to \infty} \z_{II} =  \sqrt{\frac{1}{\alpha (1-f)^2}}(\frac{\alpha}{1 + \alpha})\kappa
\end{equation}
Plugging (\ref{E:z_II_warmup}) into (\ref{E:P_base_4}) gives 
\begin{equation} 
\lim_{\V \to \infty} P(\v \text{ is lost in }[t_0,t_1]) = P(w[\Twild] = 0)
\end{equation}
We note that in Theorem \ref{T:warmup} the scaling will force the following duality in $\z_{II}$,
\begin{equation}  \label{E:warmup_z_scaling}
\lim_{\V \to \infty} \z_{II} = \bigg\{
\begin{array}{cc}
\infty & \text{if } \frac{\phi}{\beta} \ge 1. \\
0 & \text{if } \frac{\phi}{\beta} < 1.
\end{array}
\end{equation}
Since $w[\infty] = 0$ and $w[0] = 1$ we have,
\begin{equation}  \label{E:lost_prob}
P(\v \text{ is lost in }[t_0,t_1])
= \bigg\{
\begin{array}{cc}
1 & \text{if } \frac{\phi}{\beta} \ge 1. \\
0 & \text{if } \frac{\phi}{\beta} < 1.
\end{array}
\end{equation}
and this is the essential difference between the results in Theorems \ref{T:warmup} and \ref{T:main}.  If $\v$ is lost during Stage II, it is straightforward to show that the system must go to $\u_M$.  

\vspace{.1cm}
\begin{flushleft}
\textbf{Stage III:}\\
\end{flushleft}
If the wild type is not lost in Stage II,  then by the same arguments used in  Lemma \ref{L:I_and_III_s} to show that $\u \approx \bu$ in Stage I, we can show that $\u \approx \bu$ in Stage III (with the precise statement being identical to that found in (\ref{E:I_det})).

\vspace{.1cm}
\begin{flushleft}
\textbf{Stage IV:}\\
\end{flushleft}
The arguments of Stage IV are almost identical to Stage II, except that $\vs$ is now what collapses rather than $\v$.  We define a function $\psi$ in section \ref{S:stochastic} that plays the role $\phi$ did in Stage II.  And the relationship between $\psi$ and $\psilim$ is completely analogous to the relationship between $\phi$ and $\philim$.  If the mutant is lost then the system goes to $\u_W$.

\vspace{.1cm}
\begin{flushleft}
\textbf{Behavior over a full cycle:}\\
\end{flushleft}
If both the mutant and the wild type survive the first cycle of Stages I-IV then our claim is that the system gets arbitrarily close, at least to $O(1)$, to $\u_C$.   This is a result of the damped oscillations seen in Figure \ref{F:full_ODE}.  To demonstrate the damping of the oscillations and their impact on the probability of losing $\v$ or $\vs$ in a given cycle, we first note that by Lemmas \ref{L:I} and \ref{L:III}, $\bp$ changes by only $O(\epsilon |\log \epsilon|)$ during Stages I and III.  And since $\u$ and $\bu$ are linked during those stages through (\ref{E:I_det}) the same will be true of $\p$ and $\bp$.  In Stage II, since $\v$ is $O(\epsilon^q)$ we will have $\dot{\p} = -\epsilon \alpha \p^2 + O(\epsilon^{q+1})$.  In Stage IV we will show that while $\vs \le O(\epsilon^q)$ we have $\v - (1 - \p) = O(\epsilon)$.  From this we have $\dot{\p} = \epsilon \p(1 - (1 + \alpha) \p) + O(\epsilon^2)$.  Both Stages II and IV can be shown to be of duration $O(\frac{1}{\epsilon})$ and so we can explicitly integrate $\p$ through the cycle with an error term $O(\epsilon)$.  In this case we find $0 < p(\Tvsb) - (1-f) < p(\TO) - (1-f) > 0$.  Essentially we have shown that the starting point of our cycle is damped towards $(1-f)$ in subsequent cycles.  

	To connect to the probability of loss, consider the impact of this damping on  $\philim$.   $\philim$ can be thought of as a function of $\p(\TO)$ since $\phi(\frac{\se(\p(\TO)}{\epsilon})$ is, through this connection we can consider how $\philim$ changes with $\p(\TO)$.  Explicit differentiation gives,
\begin{equation}
\frac{\partial \philim(p(\TO))}{\partial \p(\TO)} = \frac{p(\TO)-(1-f)}{p^2(\TO)} < 0
\end{equation}
So with each progressive cycle, $\philim$ will be $O(1)$ smaller.  The same can be said for $\psi$.  To see the effect of this, let $\philim^{(1)}$ be the value of $\philim$ corresponding to the first cycle, i.e. exactly $\philim$, and $\philim^{(2)}$ be the value of $\philim$ corresponding to the second cycle.  Then we have,
\begin{equation}
\exp[\frac{\philim^{(2)}}{\epsilon}] = \exp[\frac{\philim^{(1)}}{\epsilon}]\bigg(\exp[-O(\frac{1}{\epsilon})]\bigg)
	= \exp[\frac{\philim^{(1)}}{\epsilon}]\bigg(\exp[-O(\log \V)]\bigg).
\end{equation}
If we follow the arguments that led to (\ref{E:lost_prob}) we see that if $\v$ is not lost in the first Stage II it will not be lost in subsequent stages, and similarly for $\vs$.  

	Now note that we set $t_f = O(\frac{1}{\epsilon^2})$.  We show in sections \ref{S:deterministic} and \ref{S:stochastic} that a full cycle takes $O(\frac{1}{\epsilon})$, so if $\v$ and $\vs$ are not lost we are considering $O(\frac{1}{\epsilon})$ cycles.  Since we have shown that $p(\Tvsb) - (1-f) < p(\TO) - (1-f)$ we see that with every cycle the next $\p(\TO)$ gets $O(1)$ closer to $1-f$.  This implies that by the time interval $[O(\frac{1}{\epsilon^\frac{3}{2}}), t_f]$, $\p - (1-f)$ will be $o(1)$.  It is then easy to show that $\v$ and $\vs$ get within $o(1)$ to their values in $\u_C$.  And so $\u$ gets within $o(1)$ of $\u_C$.  
	
	Finally, our analysis of the different stages contained expansions with error terms $o(\epsilon)$ and excluded sets with probability $O(\epsilon^2)$.  Since by time $t_f$ we have gone through at most $O(\frac{1}{\epsilon})$ cycles, if we now take $\V \to \infty$ which takes $\epsilon \to 0$ our arguments become true with probability approaching one.
\end{proof}

	Before finishing with Theorems \ref{T:warmup} and \ref{T:main} we return to an issue raised in section \ref{S:bd}.  In that section, we stated that the stochastic terms involving $\p$ had been dropped in (\ref{E:SDE}), but that our results would be unaffected.  Indeed, $\v$ and $\vs$ are only affected by stochastic terms when they drop to low levels in Stages II and IV.  $\p$ never experiences such bottlenecks, so the stochastic terms in the $\p$ equation will have no effect.  If we included them we would need a Lemma similar to Lemma \ref{L:I_and_III_s}.  But this would just add technicalities to our discussion.

	
	We now prove Theorem \ref{T:genetics}.  Before proceeding we explain how we build sample lineages on the population process $\u(t)$.  These justifications are similar to arguments found in \cite{Wakeley_2009_TPB}. Take two times $t$ and $t + \Delta t$.  Suppose at time $t + \Delta t$ we have $n'$ sample lineages of wild type.  If a birth event happens in the underlying birth-death process during $[t, t + \Delta t]$ then the number of lineages may drop from $n'$ to $n'-1$.  Indeed, if the new children formed by the birth event are both part of the $n'$ samples  at time $t + \Delta t$ then we have a coalescent event and at time $t$ there will be $n'-1$ lineages. By symmetry, given a birth event there is a probability $\frac{n'(n'-1)}{(\V \v(t))^2}$ that a coalescent will occur (recall that $\V \v(t)$ is the number of wild type cells at time $t$).

	To compute the probability of a coalescent event, we need to know the probability of a birth event in $[t, t+\Delta t]$ conditioned on the value of $\u(t + \Delta t)$.  Set $\u(t + \Delta t) = \hat{\u} = (\hat{\v}, \hat{\vs}, \hat{\p})$.  Then an application of Bayes rule gives
\begin{align}
P(\text{birth} \ | \ \u(t+\Delta t) = \hat{u}) 
	& = \frac{P(\text{birth} \ | \ \u(t) = (\hat{v} - \frac{1}{\V}, \hat{\vs}, \hat{\p})) P(\u(t) = (\hat{v} - \frac{1}{\V}, \hat{\vs}, \hat{\p}))}{P(\u(t+\Delta t) = \hat{u})}
\\ \notag
	& = \frac{\V\hat{v} - 1}{P(\u(t+\Delta t) = \hat{u} \ | \ \u(t) = (\hat{v} - \frac{1}{\V}, \hat{\vs}, \hat{\p}))}
\\ \notag
	& = \V\hat{v} - 1 + O(\Delta t).
\end{align}
	
	From these arguments we have that the probability of a coalescent event in $[t, \Delta t]$ is $\frac{n'(n'-1)}{\V \v(t)} + O(\frac{1}{(\V \v(t))^2}) + O(\Delta t)$.  We note that this approximation breaks down when $\V \v(t) = O(1)$.  However, we only use these rate computations if the wild type is not lost, otherwise wild type lineages do not exist.  And by the computations of Lemma \ref{L:stochastic_II_part2}, in this case we always have $\V \v = O(\frac{1}{\sqrt{\epsilon}})$ so our expansions are valid.  
	
	Finally, before proceeding to the proof of Theorem \ref{T:genetics} we connect our lineage distribution to the Kingman coalescent. A well-known result in coalescent theory states that if lineages coalesce at rate, say, $r(t)$ over an interval $[0,T]$, then the lineages at time zero of this time varying coalescent process will have the same distribution as the Kingman coalescent run for time $\int_0^T ds r(s)$  \cite{Tavare_2001_St_Flour_Lecture_Notes}.  Essentially, in the Kingman coalescent lineages coalesce at rate $1$, so by a time rescaling one can produce coalescent events at the correct rate, $r(t)$.  For us, the consequence of this observation is that the distribution of wild type lineages will be given by a Kingman coalescent run for time $\int_0^{t_f} ds \frac{1}{\V \v(s)}$.  Hence, in the proof of Theorem \ref{T:genetics} we focus on the quantity $\int_0^{t_f} ds \frac{1}{\V \v(s)}$

\begin{proof}[Proof of Theorem \ref{T:genetics}]
We recall the notation $(n_0, B_1, B_2, \dots, B_{n_0})$, introduced in section \ref{S:results}, representing the number of lineages left at time zero from $n$ samples at time $t_f$.  We first consider the following cases: mutant fails, wild type loss, and mutant loss.

If the mutant fails then Lemma \ref{L:I_0} shows that $\v(t) = \frac{\alpha}{1 + \alpha} + O(\epsilon^q)$ for all $t \in [0, t_f]$.  Then we have,
\begin{equation}
\lim_{\V \to \infty} \int_0^{t_f} ds \frac{1}{\V \v(s)} = \lim_{V \to \infty} O(\frac{1}{\epsilon^2 \V} = 0
\end{equation}

If the wild type is lost then, since there is only a single mutant at time $0$,  $n_0=1$ and $B_1 = n$.  This is the same distribution at $\Pi(\infty; n)$.

Now we consider the case of the mutant type being lost.  First we claim,
\begin{equation}
\lim_{\V \to \infty} \int_0^{t_f} ds \frac{1}{\V \v(s)} = 
\lim_{\V \to \infty} \int_{t_0}^{t_1} ds \frac{1}{\V \v(s)}
\end{equation}
where recall that $[t_0, t_1]$ is an subinterval in the Stage II in the first cycle of Stages I-IV.  To see this, first note that in Stages I, III, and IV and in the initial stage we have $\v \ge \epsilon^q$.  So since the duration of these stages is $O(\frac{1}{\epsilon})$ we have,
\begin{equation}
\lim_{\V \to \infty} \int_{[0,\Tvd], [\Tvb,\Tvsb]} ds \frac{1}{\V \v(s)} = \lim_{\V \to \infty} O(\frac{1}{\epsilon^{q+1} \V}) = 0.
\end{equation}
On $[\Tvsb, t_f]$ the mutant is already lost and $\v(t) = O(1)$, so $\lim_{\V \to \infty} \int_{[\Tvsb,t]} ds \frac{1}{\V \v(s)} = 0$.  We are left to consider $[\Tvd, \Tvb]$ during the first cycle.  The key relation is (\ref{E:final_vA_formula}) in Lemma \ref{L:stochastic_II_part2} which allows us to compute $\v(t)$ for $t \in [t_0, t_1]$.  By plugging in $t = \frac{\se}{\epsilon}$ in (\ref{E:final_vA_formula}) we arrive at,
\begin{equation}  \label{E:small_v_t0}
\v(t_0) = \v(\frac{\se}{\epsilon}) \exp[\alpha(1-f)^2 \frac{1}{\epsilon^{2\kval-1}}]
\end{equation}
We now wish to demonstrate the following bound,
\begin{align}  \label{E:stage_II_gen_asymp}
\int_{\Tvd}^{\Tvd} ds \frac{1}{\V \v(s)} & = \int_{t_0}^{t_1} ds \frac{1}{\V \v(s)} + O(\frac{1}{\epsilon \V \v(t_0)})
=  \int_{t_0}^{t_1} ds \frac{1}{\V \v(s)}(1 + O(\frac{\exp[\alpha(1-f)^2 \frac{1}{\epsilon^{2\kval-1}}]}{\epsilon})).
\end{align}
To see this first note that $\dot{\v}$ is always bounded by $1$.  So at least for one time unit, $\v(t)$ is of order $\v(\frac{\se}{\epsilon})$.  Then we recall from the proof of Theorem \ref{T:warmup} that $\v$ is well approximated by $\bv$ in Stage II outside of $[t_0, t_1]$, and that $\bv$ is strictly decreasing on $[\bTvd, t_0]$ and strictly increasing on $[t_1, \bTvb]$.  Since Stage II is of duration $O(\frac{1}{\epsilon})$ we must have,
\begin{equation}
\int_{[\Tvd, t_0], [\Tvb, t_1]} ds \frac{1}{\V \v(s)} \le \frac{1}{\epsilon \V \v(\frac{\se}{\epsilon})}.
\end{equation}
From the above bound and (\ref{E:small_v_t0}), (\ref{E:stage_II_gen_asymp}) follows.  Summarizing, we have shown,
\begin{equation}
\lim_{\V \to \infty} \int_0^{t_f} ds \frac{1}{\V \v(s)} = \lim_{\V \to \infty} \int_{t_0}^{t_1} ds \frac{1}{\V \v(s)}.
\end{equation}
Finally, Lemma \ref{L:stochastic_II_part2} shows $\int_{t_0}^{t_1} ds \frac{1}{\V \v(s)} \to \Upsilon \int_0^{\sqrt{2\pi}} ds \frac{1}{w_1[\Upsilon (k+1)\gamma s]}$.  

We have left the case of the wild type and mutant type both surviving.  In this case, the computations are simply a combination of the wild type loss and mutant type loss case after we split the sample into wild and mutant samples.  We have $\v(t_f) = \alpha(1-f)$ and $\vs(t_f) = f - \alpha(1-f)$.  So the probability of drawing a wild type is, as $\V \to \infty$, exactly $\frac{\alpha(1-f)}{f}$ and the number of wild types out of $n$ sample is binomial with $n$ trials and success probability $\frac{\alpha(1-f)}{f}$.
 \end{proof}

\section{The Deterministic System}  \label{S:deterministic}
\setcounter{equation}{0}
\setcounter{lemma}{0}

	In this section we consider the deterministic system (\ref{E:vvsp_system}).  Our goal will be to find asymptotic expansions for $\bu$ at the end of each of Stages I-IV and to develop estimate of the $\bu$ dynamics during the stages as well.  We define $\bdelta(t) = \bp(t) - (1-f)$.

\subsection{Stage I}

	Recall that Stage I is given by the interval $[\bTvd, \bTvb]$.  In this subsection, for notational convenience, we set $\delta = \delta(\TO)$.  We assume $\bvs(\TO) = \epsilon^q$, $\bv(\TO) - (1 - \bp(\TO)) \le O(\epsilon)$ (this assumption will connect to Lemma \ref{L:I_0} found in section \ref{S:stochastic}.

\begin{lemma}  \label{L:I}
$\bTvd - \TO = O(\frac{|\log \epsilon|}{\delta})$ and,
\begin{gather}  \label{E:I_v_vs}
|f - \bvs(\bTvd)| = O(\epsilon^q),  
\end{gather}
\begin{equation} \label{E:I_p}
\bp(\bTvd) = \bp(\TO) 
	+ \frac{q \epsilon |\log \epsilon|}{\delta} \Delta \bar{p}_{I,1}
	+ \frac{\epsilon}{\delta} \Delta \bar{p}_{I,2} + O(\epsilon^\frac{3}{2} |\log \epsilon|^2)
\end{equation}
where,
\begin{gather}
\Delta \bar{p}_{I,1} = \bp(\TO)(1 - (2 + \alpha) \bp(\TO)), \\ \notag
\Delta \bar{p}_{I,2} = \bp(\TO)(1 - (1+\alpha)\bp(\TO)) \log(f) - \alpha \bp^2(\TO) \log(1 - \bp(\TO)).
\end{gather}
\end{lemma}
\begin{proof}
We will separate $[\TO, \bTvd]$ into three intervals using stopping times $T_1, T_2$ where,
\begin{gather}
T_1 = \inf\{t-\TO : \bvs(t) = \sqrt{\epsilon}\}. \\
T_2 = \inf\{t-\TO : \bv(t) = \sqrt{\epsilon}\}
\end{gather}
Our first task is to determine $T_1$.  Consider $t \in [\TO,\bTvd]$. For such $t$ we have,
\begin{equation}  \label{E:bound_1_p_v}
\dot{(1 - \bp(\TO) - \bv)} = -\bv(1 - \bp - \bv - \bvs) 
	= -\bv(1 - \bp(0) - \bv - \bvs - O(\epsilon t)).
\end{equation}
where we have used the fact that $\bp(t) - \bp(\TO) \le \epsilon t$.  From (\ref{E:bound_1_p_v}) we see that if $(1 - \bp(\TO) - \bv) > O(\max(\epsilon t, \bvs, \epsilon))$ then  $\dot{(1 - \bp(\TO) - \bv)} < 0$.  Since  $(1 - \bp(\TO) - \bv(\TO)) = O(\epsilon)$ we can conclude,
 $|1 - \bp(\TO) - \bv| = O(\max(\epsilon t, \bvs,\epsilon))$.  Using (\ref{E:bound_1_p_v}), we have for $t \in [\TO, \TO + (T_1 \wedge O(|\log \epsilon|))]$
\begin{gather}  \label{E:early_dbvs}
\dbvs = \dbvs(f - \bv - \bvs) \ge \dbvs(\delta + O(\sqrt{\epsilon}) 
\end{gather}
From (\ref{E:early_dbvs}), we see that $T_1  = O(|\log \epsilon|)$ and so we have $\exp[(\delta + O(\sqrt{\epsilon}))T_1] = \frac{\sqrt{\epsilon}}{\epsilon^q}$.  Solving for $\T_1$ we arrive at $T_1 = \frac{q-\frac{1}{2}}{\delta}|\log \epsilon|  + O(\sqrt{\epsilon} |\log \epsilon|^2)$. Integrating $\bp$ we find,
\begin{align}
\bp(T_1) & = \bp(\TO) + \epsilon \bp(\TO)(1-(1 + \alpha) \bp(\TO))T_1 + O(\epsilon^\frac{3}{2} |\log \epsilon|^3)
\\ \notag
  & = \bp(\TO) + \frac{(q-\frac{1}{2}) \epsilon |\log \epsilon|}{\delta} \bp(\TO)(1-(1 + \alpha)\bp(\TO)) +  O(\epsilon^\frac{3}{2} |\log \epsilon|^3)
\end{align}

	Now we consider the time interval $[\TO + T_1,\TO + T_2]$.  Let $r(t) = \frac{(1 - \bp(\TO)) - \bv(t)}{\bvs(t)}$.  Then direct computation gives,
\begin{equation}
\dot{r} = (1 - \bp(\TO)) - fr + \frac{\bv}{\bvs} (\bp - \bp(\TO)).
\end{equation}
Integrating through we have,
\begin{align}  \label{E:integrated_r}
r(t) = & r(\TO) \exp[-ft] + (1-\bp(\TO))\int_{\TO}^t ds \exp[-f(t-s)]
\\ \notag
	&  + \int_{\TO}^t ds \exp[-f(t-s)](\bp(s) - \bp(\TO)) (\frac{\bv}{\bvs})(s)
\end{align}
By differentiating $\log \bv + \log \bvs$ we can arrive at,
\begin{equation}
\frac{\bv}{\bvs}(t) = \frac{\bv}{\bvs}(\TO) \exp[-\delta (t-\TO) + O(\epsilon t)]
\end{equation}
Plugging this relation into the second integral term in (\ref{E:integrated_r}) gives
\begin{equation}  \label{E:integral_bound}
|\int_{\TO}^t ds \exp[-f(t-s)](\bp(s) - \bp(\TO)) (\frac{\bv}{\bvs})(s)|
	\le O(1) \exp[-\delta (t-\TO)] \frac{\sup_{s \in [\TO,t]} |\bp(s) - \bp(\TO)|}{\bvs(\TO)}
\end{equation}
Noting that $r(\TO) = O(1)$, $\bp(t) - \bp(\TO) \le \epsilon t$ and using (\ref{E:integral_bound}) in (\ref{E:integrated_r}) with $t = \TO + T_1$ gives,
\begin{align}
r(\TO + T_1) & = \frac{1 - \bp(\TO)}{f} + O(\epsilon^\frac{f(q-\frac{1}{2})}{\delta}) 
\\ \notag
	& = \frac{1 - \bp(\TO)}{f} + O(\sqrt{\epsilon} |\log \epsilon|^3).
\end{align}
Now again we use (\ref{E:integrated_r}), but instead of integrating from $0$ to $t$ we integrate from $\TO + T_1$ to $t \in (\TO + T_1, \TO + T_2)$.  In this case, since $\bvs \ge \sqrt{\epsilon}$,
\begin{equation}  
|\int_{\TO}^t ds \exp[-f(t-s)](\bp(s) - \bp(\TO)) (\frac{\bv}{\bvs})(s)|
	\le O(\sqrt{\epsilon}|\log \epsilon|).
\end{equation}
which leads to
\begin{equation}  \label{E:r_bound_for_bvs}
r(t) = \frac{1 - \bp(\TO)}{f} + O(\sqrt{\epsilon} |\log \epsilon|^3). 
\end{equation}
Exploiting (\ref{E:r_bound_for_bvs}) allows us to integrate $\bv$ explicitly.  We have,
\begin{gather}
\bv(\TO + T_1) = (1-\bp(\TO)) - (\frac{1 - \bp(\TO)}{f})\sqrt{\epsilon} + O(\epsilon |\log \epsilon|^3). \\ 
\bvs(t) = f - (\frac{f}{1-\bp(\TO)}) \bv(t) + O(\epsilon |\log \epsilon|^3).
\end{gather}
Substituting these expressions into (\ref{E:vvsp_system}) gives,
\begin{equation}  \label{E:logistic_v}
\dbv = -\frac{\delta}{1 - \bp(\TO)} \bv \bigg((1 - \bp(\TO) + O(\epsilon |\log \epsilon|^3)) - \bv\bigg).
\end{equation}
We can explicitly integrate (\ref{E:logistic_v}).
\begin{equation}
\bv(t) = \frac{((1-\bp(\TO))^2 + O(\sqrt{\epsilon}))\exp[-\delta(t-T_1)]}
				{\sqrt{\epsilon}(\frac{1-\bp(\TO)}{f}) + O(\epsilon |\log \epsilon|^2)
					+ (1-\bp(\TO) - \sqrt{\epsilon}(\frac{1-\bp(\TO)}{f}) + O(\epsilon |\log \epsilon|^2))\exp[-\delta(t-T_1)]}
\end{equation}
Setting $\bv(T_2) = \sqrt{\epsilon}$ and solving for $t$ gives,
\begin{equation}
T_2 - T_1 = \frac{1}{\delta} |\log \epsilon| - \frac{1}{\delta} \log(f(1-\bp(\TO))) + O(\sqrt{\epsilon}).
\end{equation}
and explicit integration gives,
\begin{equation}  \label{E:int_v_T_1_T_2}
\int_{T_1}^{T_2} ds \bv(s) 
	= \frac{1-\bp(\TO)}{2\delta} |\log \epsilon| + \frac{1-\bp(\TO)}{\delta} \log (f) + O(\sqrt{\epsilon})
\end{equation}
Using (\ref{E:int_v_T_1_T_2}) we find,
\begin{align}
\bp(\TO + T_2) = & \bp(\TO + T_1) + \frac{\epsilon |\log \epsilon|}{\delta} \bigg(\frac{\bp(\TO)(1-\bp(\TO))}{2} 
						- \alpha \bp^2(\TO)\bigg)
\\ \notag
	& + \frac{\epsilon}{\delta} \bigg(\bp(\TO)(1-\bp(\TO)) \log (f) - \alpha \bp^2(\TO) \log(f(1-\bp(\TO)))\bigg) 
			 + O(\epsilon^\frac{3}{2})
\end{align}

	Finally we consider $[\TO + T_2,\bTvd]$. The analysis in this case is almost identical to that of $[\TO,\TO + T_1]$ and we find:
$\bTvd - T_2 = \frac{q-\frac{1}{2}}{\delta}|\log \epsilon| + O(\sqrt{\epsilon})$ and 
\begin{equation}
p(\bTvd) = \bp(\TO + T_2) - \frac{\epsilon |\log \epsilon|}{\delta} \bigg(\alpha \bp^2(\TO) (q-\frac{1}{2})\bigg) + O(\epsilon^\frac{3}{2} |\log \epsilon|^3).
\end{equation}

Combining our estimates of $\bp(\TO + T_1), \bp(\TO + T_2), \bp(\bTvd)$ gives (\ref{E:I_p}).  Finally we note that our arguments show that for all $t \in [\TO, \Tvd]$ we have $\bv(t), \bvs(t) \ge \epsilon^q$.  
 \end{proof}

\subsection{Stage II}

From Lemma \ref{L:I} we can assume $\delta(\bTvd) > 0$, $\bv(\bTvd) = \epsilon^q$, and $|\bvs(\bTvd) - f| = O(\epsilon^q)$.  

\begin{lemma}  \label{L:II}
Let $\HII$ be the solution of $(1-f)\HII = \frac{1}{\alpha} \log(1 + \alpha \p(\bTvd) \HII)$.   Then,
\begin{equation}
|(\Tvb - \Tvd) - \frac{\HII}{\epsilon}| = O(\epsilon^{q-1})
\end{equation} 
\begin{equation}
\bp(\Tvb) = \frac{\bp(\Tvd)}{1 + \alpha \bp(\Tvd) \HII} + O(\epsilon^{q}).
\end{equation}
For $t \in [\Tvd, \Tvb]$ we have,
\begin{gather}
\sup_{t} |\bvs(t) - f| \le O(\epsilon^q) \\ \notag
\bv(t) = \bv(\Tvd) \exp\bigg[(1 - f)(t - \Tvd) - \frac{1}{\alpha \epsilon}\log\big(1 + \epsilon \alpha \bp(\Tvd) (t - \Tvd)\big) + O(\epsilon^{q-1})\bigg]
\end{gather}
\end{lemma}

\begin{proof}
We first show that on $[\bTvd, \bTvb]$, $\bp$ has simple dynamics.  Indeed, for all $t \in [\bTvd, \bTvb]$, 
\begin{equation}
\ddbp = -\epsilon \alpha \bp^2 + C \epsilon^{1 + q},
\end{equation}
From this we have,
\begin{equation}  \label{E:p_equation}
\bp(t) = \frac{\bp(\Tvd)}{1 + \epsilon \alpha \bp(\Tvd) (t - \Tvd)} + O(\epsilon^q)
\end{equation}
Now consider $\bvs$ for $t \in [\Tvd, \Tvb]$.  
\begin{equation}  \label{E:bvs_top}
\dot{f - \bar{\v}^*} = -\bvs(f - \bv - \bvs) = -\bvs(f - \bvs - C \epsilon^{q}).
\end{equation}
where $C \le 1$.  Since $|f - \bvs(\Tvd)| = O(\epsilon^q)$, from (\ref{E:bvs_top}) we can conclude that  $|f - \vs| \le O(\epsilon^q)$ on $[\Tvd, \Tvb]$.

	Now we bound $\Tvb$.  We have for $t \in [\Tvd, \Tvb]$,
\begin{equation}
\dbv = \bv(1 - f - \bp + O(\epsilon^q))
\end{equation}
Plugging (\ref{E:p_equation}) into the relation directly above leads to
\begin{equation}
\bv(t) = \v(\Tvd) \exp[(1-f)(t - \Tvd) - \frac{1}{\alpha \epsilon} \log(1 + \epsilon \alpha \bp(\Tvd) (t - \Tvd)) + O(\epsilon^q t)]
\end{equation}
Since $\bp(\Tvd) > 1-f$, a simple Taylor series argument shows,
\begin{equation}
v(\Tvd + \frac{\HII}{\epsilon} - O(\epsilon^{q-1})) < \epsilon^q < v(\Tvd + \frac{\HII}{\epsilon} + O(\epsilon^{q-1})).  
\end{equation}
So we can conclude that $\Tvd + \frac{\HII}{\epsilon} - O(\epsilon^{q-1}) < \Tvb < \Tvd + \frac{\HII}{\epsilon} + O(\epsilon^{q-1})$ and the lemma follows.
 \end{proof}

\subsection{Stage III}

	From Lemma \ref{L:II} we can assume $\delta(\bTvb) < 0$, $\bv(\bTvb) = \epsilon^q$, and $|\bvs(\bTvb) - f| \le O(\epsilon^q)$.  For convenience in this section set $\delta = \delta(\bTvb)$.  

\begin{lemma}  \label{L:III}
$\bTvsd - \bTvb = O(\frac{|\log \epsilon|}{|\delta|})$,
\begin{gather}  \label{E:III_v_vs}
|\bv(\Tvsd) - (1 - \bp(\Tvsd))| = O(\epsilon),  
\end{gather}
\begin{equation} \label{E:III_p}
\bp(\bTvsd) = \bp(\bTvb)
	+ \frac{q \epsilon |\log \epsilon|}{|\delta|} \Delta \bar{p}_{III,1} 
	+ \frac{\epsilon}{|\delta|} \Delta \bar{p}_{III,2} + O(\epsilon^\frac{3}{2} |\log \epsilon|^2)
\end{equation}
where,
\begin{gather}
\Delta \bar{p}_{III,1} = \bp(\bTvb)(1 - (2 + \alpha) \bp(\bTvb)), \\ \notag
\Delta \bar{p}_{III,2} = \bp(\bTvb)(1 - (1+\alpha)\bp(\bTvb)) \log(f) - \alpha \bp^2(\bTvb) \log(1 - \bp(\bTvb)).
\end{gather}
\end{lemma}

\begin{proof}
The proof of the following lemma is almost identical to the of Lemma \ref{L:I}.  The only difference involves the proof of (\ref{E:III_v_vs}).  Set $\Delta \bv = \frac{\bv - (1-\bp)}{\epsilon}$.  Tnen direct computation gives,
\begin{equation}  \label{E:delta_v_equation}
\dot{\Delta \bv} = -\bv \Delta \bv - \frac{\bv \bvs}{\epsilon} + \bp(\bv - \alpha \bp)
\end{equation}
We define $T_1, T_2$ in the same manner as in Lemma \ref{L:I}, except that the role of $\bv$ and $\bvs$ are interchanged.  That is, $\bv(T_1) = \sqrt{\epsilon}$, $\bvs(T_2) = \sqrt{\epsilon}$.  Essentially the same arguments we used in Lemma \ref{L:I} involving $r(t)$ can be used here to show $\Delta \bv(T_2) \le O(\frac{1}{\sqrt{\epsilon}})$.  

	We now integrate (\ref{E:delta_v_equation}) from $T_2$ to $\bTvsd$.  First notice that in this range we have,
\begin{equation}  \label{E:delta_v_equation_2}
\dot{\Delta \bv} = -\bv \Delta \bv - O(\frac{1}{\sqrt{\epsilon}}).
\end{equation}
From which we see that for $t \in [T_2, \bTvsd]$, $\Delta \bv \le O(\frac{1}{\sqrt{\epsilon}})$ and so $\bv = 1 - \bp(\bTvb) + O(\sqrt{\epsilon})$.
Now we proceed to integrate (\ref{E:delta_v_equation}) using the integrating factor $\exp[\int_{\bTvb}^t ds \bv(s)]$.  Since everything we do below is in terms of orders, we can replace this integration factor by $\exp[(1 - \bp(\bTvb))(t - T_2)]$.  This gives,
\begin{align}  \label{E:integrated_delta_v}
\Delta \bv(\bTvsd) = & O(\Delta \bv(T_2) \exp[-(1 - \bp(\bTvb))(\bTvsd - T_2)])
\\ \notag
		& + \frac{1}{\epsilon} \int_{T_2}^{\bTvsd} ds \bvs(s) \exp[(1 - \bp(\bTvb))(\bTvsd - s)]
		+ O(1).
\end{align}
By the same arguments as in Lemma \ref{L:I} we can show $\bvs(s) = O(\sqrt{\epsilon}\exp[-\delta (s - T_2)])$ and $\bTvsd - T_2 = \frac{(q - \frac{1}{2})}{\delta}|\log \epsilon|$.  Plugging all this into (\ref{E:integrated_delta_v}) finally leads to $\Delta v(\bTvsd) = O(1)$.
 \end{proof}

\subsection{Stage IV}

From the results of Lemma \ref{L:III} we can assume $\delta(\bTvsd) < 0$, $|\bv(\bTvsd) - (1 - p(\bTvsd))| = O(\epsilon)$, and $\bvs = \epsilon^q$.  

\begin{lemma}  \label{L:IV}
Define $\HIV$ as the solution of the following equation,
\begin{equation}
(1-f)\HIV = \frac{1}{1 + \alpha} \log\bigg(1 + (1 + \alpha)\bp(\bTvsd)(\exp[\HIV]-1)\bigg).
\end{equation}
Then $|(\bTvsb - \bTvsd) - \frac{\HIV}{\epsilon}| = O(\epsilon)$ and
\begin{equation}
\bp(\bTvsb) = \frac{\bp(\bTvsd) \exp[\HIV]}{1 + (1 + \alpha) \bp(\bTvsd)(\exp[\HIV]-1)} + O(\epsilon).
\end{equation}
For $t \in [\bTvsd, \bTvsb]$ we have,
\begin{gather}
|\bv(t) - (1 - \bp(t)| \le O(\epsilon) \\ \notag
\bvs(t) = \bvs(\bTvsd) \exp[g(t) + O(\epsilon^2 t)]
\end{gather}
where,
\begin{align}
g(t) = & -(1 - f)(t - \bTvsd) + \frac{1}{(1 + \alpha)\epsilon}
	\log\bigg(1 + (1 + \alpha)\bp(\bTvsd)(\exp[\epsilon (t-\bTvsd)]-1)\bigg)
\\ \notag
	& + \epsilon (\frac{1 - \bp(\bTvsd)}{1 - \bp(t)}) \exp[\bp(\bTvsd) - \bp(t)].
\end{align}
\end{lemma}

\begin{proof}
We start by considering $\bv - (1 - \bp)$.  Recall the definition of $\Delta \bv$ from Lemma \ref{L:III}.  By the same arguments as in Lemma \ref{L:III} we can show that $\Delta \bv \le O(1)$ and hence $\bv = (1 - \bp) + O(\epsilon)$.  However, to compute $\bvs$ to $O(\epsilon)$ we need $\bv$ to one higher order.  To see why we need $\bv$ to one higher order, suppose we have $\bv = 1 - \p + O(\epsilon)$.  Then we would have the estimate $\dot{\bvs} = \bvs(\p - (1-f) + O(\epsilon)$ and integrating this expression for $O(\frac{1}{\epsilon})$ time units leads to $O(1)$ terms.   To estimate $\bv$ to $O(\epsilon)$ we integrate (\ref{E:delta_v_equation}) and apply our apriori bound of $\Delta \bv = O(1)$ to arrive at,
\begin{align}
\Delta \bv(t) = & (1 - \bp(t)) + I_1 + I_2 + O(\epsilon^2)(t - \bTvsd),
\end{align}
where,
\begin{gather}
I_1 =  \epsilon \Delta \bv(0) \exp[-\int_{\bTvsd}^t ds(1 - \bp)] \\ \notag
I_2 = \epsilon \int_{\bTvsd}^t ds \exp[-\int_{s}^t ds' (1 - \bp)] p(1 - (1+\alpha)p) \notag
\end{gather}
$I_1$ is $O(\epsilon)$ for all $t \in [\bTvsd, \bTvsb]$ since $1 - \bp(t) > O(1)$.  Consider $I_2$.  We have $\dot{\bp} = \epsilon \bp(1 - (1+\alpha)\bp) + O(\epsilon^2)$.  So we have,
\begin{equation}
I_2 = \int_{\bTvsd}^t ds \exp[-\int_{s}^t ds' (1 - \bp)] \dot{\bp}
\end{equation}
We now apply an integration by parts trick to obtain an expansion for $I_2$.  Indeed, since every derivative of $\bp$ earns us another $\epsilon$, integration by parts gives the following.
\begin{align}
I_2 & = \int_{\bTvsd}^t ds (1 - \bp) \exp[-\int_{s}^t ds' (1 - \bp)] \frac{\dot{\bp}}{1 - \bp}
\\ \notag
	& = \bigg(\exp[-\int_{s}^t ds' (1 - \bp)] 
	\frac{\dot{\bp}(s)}{1 - \bp(s)}\bigg)\bigg|_{s=\Tvsd}^{s=t} + O(\epsilon^2)(t - \bTvsd)
\\ \notag
	& = \frac{\dot{\bp}(t)}{1 - \bp(t)} + O(\epsilon^2)(t - \bTvsd)
\end{align}
So we have the expansion
\begin{equation}
\bvs = (1 - \bp) + \epsilon \frac{\bp(t)(1 - (1+\alpha)\bp(t))}{1 - \bp(t)} + O(\epsilon^2)(t - \bTvsd).
\end{equation}

We now consider $\bp$.  First note that in Stage IV we have,
\begin{equation}
\dot{\bp} = \epsilon \bp(1 - (1+\alpha)\bp) + O(\epsilon^2).
\end{equation}
Solving for $\bp$ gives,
\begin{equation}  \label{E:p_up_formula}
\bp(t) = \frac{\bp(\bTvsd) \exp[\epsilon (t-\bTvsd)]}{(1 + (1 + \alpha) \p(\bTvsd) (\exp[\epsilon (t-\bTvsd)]-1)} + O(\epsilon^2)(t - \bTvsd)
\end{equation}

Now we turn to $\bvs$.  From the arguments above, by plugging our expansion for $\bv$ into $\dot{\bar{\vs}} = \bvs(f - \bv - \bvs)$, we have
\begin{equation}
\bvs(t) = \bvs(\bTvsd)\exp[g(t) + O(\epsilon)].
\end{equation}
Recall that $\bTvsb$ is defined by $\bvs(\bTvsb) = \epsilon^q$.  Considering $g(t)$ and the definition of $\HIV$, a Taylor expansion gives,
\begin{equation}
\bTvsb - \bTvsd = \frac{\HIV}{\epsilon} + O(\epsilon)
\end{equation}
and,
\begin{equation}
\bvs(\bTvsb) = \bvs(\bTvsd)\exp[g(\bTvsb) + O(\epsilon)].
\end{equation}
 \end{proof}

\section{The Stochastic System}  \label{S:stochastic}
\setcounter{equation}{0}
\setcounter{lemma}{0}

	In this section we consider the stochastic system (\ref{E:SDE}).  Similarly to section \ref{S:deterministic}, our goal will be to find asymptotic expansions for $\u$ at the end of each of Stages I-IV and develop estimates for the dynamics of $\u$ within the stages.  In this section we also consider the initial stage.    We recall from section \ref{S:deterministic} that $\delta(t) = \bp(t) - (1-f)$ and we emphasize that $\delta$ depends on the deterministic system (\ref{E:vvsp_system}) and not the stochastic system (\ref{E:SDE}).  Finally, in this section we will assume that $\epsilon$ obeys the scaling of Theorem \ref{T:main}.  The arguments in the case of Theorem \ref{T:warmup} scaling are similar, and in fact simpler.  

\subsection{Initial Stage}

Recall that at $t=0$ we have,
\begin{equation}
\v(0) = \frac{\alpha}{1 + \alpha}, \quad
\p(0) = \frac{1}{1 + \alpha}, \quad
\vs(0) = \frac{1}{\V}.
\end{equation}
Set $\Delta \v = \v - \v(0)$ and $\Delta \vs = \vs - \vs(0)$.   Then define T by 
\begin{equation}
T = \inf\{t : \vs \notin (0,\epsilon^q), |\Delta \v| = \epsilon^q,\text{or } |\Delta \p| = \epsilon^q\}
\end{equation}
Note that $T$ is not quite $\TO$, but the lemma below will show that the two are equivalent in the $\V \to \infty$ limit.

We have the following lemma.
\begin{lemma}   \label{L:I_0}
\begin{equation}
\lim_{\V \to \infty}P(\vs(T) = \epsilon^q \text{ or } \vs(T) = 0) = 1.
\end{equation}
\begin{equation}  \label{E:I_O_P}
\lim_{\V \to \infty} P(\vs(T)  = 0) = \exp[-\frac{4(f - \alpha(1-f))}{(k^*+1)(1 + \alpha)}].
\end{equation}
\begin{equation}  \label{E:I_O_T}
\lim_{\V \to \infty} P(T \le \frac{2\log \V}{f - \v(0)}) = 1.
\end{equation}
\end{lemma}

\begin{proof}
Define two stopping times,
\begin{gather}
T_\Delta = \inf\{t : |\Delta \v| = \epsilon^q \text{ or } |\Delta \p| = \epsilon^q\} \\ \notag
T_* = \inf\{t : \vs \notin (0,\epsilon^q)\} \notag
\end{gather}
Then clearly $T = \min\{T_\Delta, T_*\}$.  We now consider the following stochastic system.  Essentially this system is (\ref{E:SDE}) except that $\vs,\v$ and $\p$ becomes fixed once $\vs$ exits $(0, \epsilon^q)$, $|\Delta \v| \ge \epsilon^q$, and $|\Delta \v| \ge \epsilon^q$ respectively.  $\chi$ is the indicator function.
\begin{align}  \label{E:tilde_SDE}
d\tv = &  \chi(|\tv - \v(0)| \le \epsilon^q) \tv(1 - (\tv + \tvs) - \tp) dt \\ \notag
	& + \chi(|\tv - \v(0)| \le \epsilon^q) \sqrt{\frac{\tv(k + (\tv + \tvs) + \tp)}{\V}} dB_1(t),
\end{align}
\begin{equation} \notag
d\tvs = \chi(\tvs \le \epsilon^q) \tvs(f - (\tv + \tvs)) dt
	+ \chi(\tvs \le \epsilon^q) \sqrt{\frac{\tvs(k^* + (\tv + \tvs))}{\V}} dB_2(t), 
\end{equation}
\begin{equation} \notag
d\tp = \chi(|\tp - \p(0)| \le \epsilon^q)\epsilon \tp(\tv - \tp) dt
\end{equation}
Notice that in the system (\ref{E:tilde_SDE}) we are guaranteed $|\tv - \v(0)| \le \epsilon^q$, $|\tp - \p(0)| \le \epsilon^q$, and $\tvs \le \epsilon^q$ for all time.  Define $\tilde{T}_\Delta$ and $\tilde{T}_*$ analogously to $T_\Delta, T_*$ and notice that $T = \min\{\tilde{T}_\Delta, \tilde{T}_*\}$ since up to time $T$ the systems (\ref{E:SDE}) and (\ref{E:tilde_SDE}) evolve identically.  

	We now consider $\tilde{T}_*$. To control $\tvs$ we bound it from above and below by two diffusions that are simpler to analyze.  Set,
\begin{gather}  \label{E:A_B}
d\tvs_A(t) = \tvs_A(f - \v(0)) dt
	+ \sqrt{\frac{\tvs_A(k^* + \tvs_A + \tp + \tv))}{\V}} dB_2(t), \\ \notag
d\tvs_B(t) = \tvs_B(f - v(0) - \epsilon^q) dt
	+ \sqrt{\frac{\tvs_A(k^* + \tvs_B + \tp + \tv))}{\V}} dB_2(t),
\end{gather}
and notice that the drift term in $\tvs$ is less than and greater than the drift terms of $\tvs_A$ and $\tvs_B$ respectively.  By pathwise uniqueness of solution we have $\tvs_B \le \tvs \le \tvs_A$ pointwise for all $t \le \tilde{T}_*$ on the probability space of the two Brownian motions $B_1,B_2$ \cite{Karatzas_and_Shreve}.  For any $t$ we have the following bound,
\begin{equation}  \label{E:sandwich_bound}
P(\tilde{T}_* < t) \ge P(\tvs_A(t) = 0) + P(\tvs_B(t) > \epsilon^q).
\end{equation}
We would like to run $\tvs_A$ and $\tvs_B$ long enough to ensure that either these diffusions will be absorbed or reach $\epsilon^q$.  If we ignore the variance terms in (\ref{E:A_B}) we would find that $\tvs_A(t) = \epsilon^q$ for $t \approx \frac{1}{f - \v(0)} \log \V$.  With this in mind, we set $\hht = \frac{2}{f - \v(0)} \log \V$ and proceed to examine $\tvs_A(\hht)$ and $\tvs_B(\hht)$.  We consider first $P(\tvs_A(\hat{t}) = 0)$.  Set $x_A = \V \exp[-(f - \v(0))t] \tvs_A$, then
\begin{equation}
dx_A = \sqrt{x_A(k^* + v(0) + O(\epsilon^q)) \exp[-(f-\v(0))t]} dB_2(t).
\end{equation}
We now perform a time change, $w_A(\tau(t)) = x_A(t)$ with 
\begin{equation}
\tau_A'(t) = (k^* + 1 + O(\epsilon^q)) \exp[-(f-\v(0))t].
\end{equation}
This leads to
\begin{equation}  \label{E:w_A}
dw = \sqrt{w}dB_2,
\end{equation}
with $w(0) = 1$.  Note that $w$ is a weak solution to (\ref{E:w_A}) since $w$ is not $B_2$ measurable.  However, we only care about the distribution of $w(\hht)$, and by weak uniqueness $w(\hht)$ will have the same distribution as the strong solution of (\ref{E:w_A}) \cite{Karatzas_and_Shreve}.  

Recall we are interested in $P(\tvs_A(\hht) = 0)$.  Let $t_A = \tau_A(\hht)$.  Then,
\begin{equation}
P(\tvs_A(\hht) = 0) = P(x_A(\hht) = 0) = P(w[t_A] = 0).
\end{equation}
Explicit integration gives
\begin{equation}
t_A  = \frac{1}{2} (\frac{k^* + 1}{f - \v(0)}) + O(\frac{\epsilon^q}{f - \v(0)})
\end{equation}
Then standard results, see \cite{Athreya_Ney_Book} p. 260, give 
\begin{equation}  \label{E:A_sandwich_bound}
P(\tvs_A(\hht) = 0) = P(w[t_A] = 0) = \exp[-\frac{2}{t_A}] = \exp[-\frac{4(f - \v(0))}{k^* + 1}] + O(\epsilon^q).
\end{equation}
Now we consider $P(\tvs_B(\hht) \ge \epsilon^q)$.  We perform the same series of transforms as we did for $\tvs_B$ except that now we have,
\begin{gather}
x_B = \V \exp[-(f - \v(0) - 2\epsilon^q)t] \tvs_B, \\ \notag
\tau_B'(t) = (\frac{2}{k^*+1+O(\epsilon^q)}) \exp[(f - \v(0) - 2\epsilon^q)t].
\end{gather}
Setting $t_B = \tau_B(\hht)$ we can compute,
\begin{equation}
t_B = t_A(1 + O(\frac{\epsilon^q}{f-\v(0)})).
\end{equation}
We can follow our transforms to find the following,
\begin{align}
P(\tvs_B(\hht) \ge \epsilon^q)  = P(w[t_B] \ge \alpha)
\end{align}
where $\alpha = \frac{\epsilon^q \V}{\exp[(f - \v(0) - O(\epsilon^q))\hht]}$.  Plugging in our scaling for $\epsilon$ gives $\alpha = O(\frac{1}{\V (\log \V)^q})$.    All this leads to,
\begin{equation}
\lim_{\V \to \infty} P(w[t_B] \ge \alpha) = \lim_{\V \to \infty} P(w[t_A] \ne 0),
\end{equation}
If we now examine the dynamics of $\tilde{\v}$ and $\tilde{\p}$, the same type of methods that we applied to $\tilde{\vs}$ show that $\T_\Delta > O(\log \V)$ with probability approaching one and (\ref{E:I_O_P}) follows.  To see (\ref{E:I_O_T}) notice that,
\begin{equation}
\lim_{\V \to \infty} P(T < \hht) = \lim_{\V \to \infty} P(w[t_A] = 0) + \lim_{\V \to \infty} P(w[t_B] \ge \alpha) = 1
\end{equation}
 \end{proof}


\subsection{Stage I and III}

 We assume that $\bu(\TO) = \u(\TO)$ and that $\bu(\TO)$ satisfies the conclusion of Lemma \ref{L:I_0}.  Then we have the following result which covers stage I behavior.  An identical analogous result exists for stage III behavior. 
 
\begin{lemma}  \label{L:I_and_III_s}
\begin{equation}
P(\sup_{\TO \le t \le \Tvd} \|u(t) - \bar{u}(t)\|_\infty \ge \frac{1}{\V^\frac{1}{8}})
	\le O(\frac{(\log \log \V)^2 (\log \V)^2}{\V^\frac{1}{4}})
\end{equation}
\end{lemma}

\begin{proof}
Set $D = \frac{1}{\V^\frac{1}{4}}$ and define
\begin{equation}
\w = \frac{\u - \bu}{D}
\end{equation}
We have,
\begin{equation}
d\w = (\frac{a(u) - a(\bar{u})}{\delta}) dt + \frac{\sigma(u)}{\sqrt{\V}D} dB
\end{equation}
where
\begin{equation}  \label{E:def_a}
a(\bu) = \bigg(
\begin{array}{c}
\bv(1 - \bv - \bvs - \bp) \\
\vs(f - \bv - \bvs) \\
\epsilon \bp (\bv - \bp)
\end{array}
\bigg), 
\sigma(\bu) = \bigg(
\begin{array}{c}
\sqrt{\bv(k + (\bv + \bvs) + \bp)} \\
\sqrt{\bvs(k^* + (\bv + \bvs))}  \\
\end{array}
\bigg)
\end{equation}
Let $T$ be a stopping time defined as $T = \min\{t > \TO :\|\w(t)\| > 1\}$ and  $\w^T(t) = \w(t \wedge T)$.  We will show that $T \gg |\log \epsilon|$ with high probability and this will lead to $T \gg \Tvd$ with high probability.    Employing a Taylor expansion, we have for all $\TO \le t < T$,
\begin{equation}
d\w^T(t) = (\nabla a(\bu)\w^T + O(D))dt + \frac{\sigma(\bu + D \w^T)}{\sqrt{V}D} dB
\end{equation}
Using the integrating factor $\exp[-\int_0^t ds \nabla a(\bu)]$ we can arrive at
\begin{align}  \label{E:w_equation_I}
\w^T(t) = & \int_0^{t \wedge T} ds \exp[\int_s^t ds' \nabla a(\bu(s'))] O(D)
\\ \notag
	&
		+ \int_0^{t \wedge T} dB(s) \exp[\int_s^t ds' \nabla a(\bu(s'))] \frac{\sigma(\bu + D \w^T)}{\sqrt{V} D}
\end{align}
Simple computaiton gives 
\begin{equation}
\nabla a(\bu) =
\bigg(
\begin{array}{ccc}
(1 - \bv - \bvs - \bp) - \bv & -\bv & -\bv \\
-\vs & (f - \bv -\bvs) - \bvs & 0 \\
\epsilon \bp & 0 & -2\epsilon \bp
\end{array}
\bigg)
\end{equation}
Since $\bv,\bvs,\bp$ are bounded we have $\| \nabla a(\bu) \| \le 1$.  Using this bound in (\ref{E:w_equation_I}) along with a standard martingale argument gives,
\begin{equation}
E[\sup_{t \le C |\log \epsilon|} \w^2(t \wedge T)] \le O(D^2 (|\log \epsilon|)^2 \exp[2  |\log \epsilon|]) 
		+ |\log \epsilon| \frac{\exp[2 |\log \epsilon|]}{\V D^2}
\end{equation}
Recalling that $\epsilon = O(\frac{1}{\log \V})$ gives,
\begin{align}  \label{E:E_bound_I}
E[\sup_{t \le C |\log \epsilon|} \w^2(t \wedge T)]
	& = O(\frac{(\log \log \V)^2 (\log \V)^2}{\sqrt{\V}}).
\end{align}
A simple Chebyshev argument now gives,
\begin{equation}  \label{E:I_stoch_bound}
P(\sup_{t \le C |\log \epsilon|} \|(\u - \bu)(t \wedge T)\|_\infty \ge \frac{1}{\V^\frac{1}{8}}) 
	\le O(\frac{(\log \log \V)^2(\log \V)^2}{\V^\frac{1}{4}})
\end{equation}
But then by the definition of $T$ we may replace $t \wedge T$ by $t$ in the expression directly above.  Finally we note that $\bTvd = O(|\log \epsilon|)$ and $\bu(t) = O(\epsilon^{q+1})$ in time $O(|\log \epsilon|)$.  Using (\ref{E:I_stoch_bound}) we see that $\Tvd = O(|\log \epsilon|)$ and so we may replace the term $\sup_{t \le C |\log \epsilon|}$ in (\ref{E:I_stoch_bound}) by $\sup_{t \le \Tvd}$.
 \end{proof}

\subsection{Stage II}

Now we consider the time interval $[\Tvd, \Tvb]$.  This interval is where stochastic effects become important and $\v$ diverges from $\bv$.  Note that from Lemmas \ref{L:I} and \ref{L:I_and_III_s} we have $|\u - \bu| \le O(\frac{1}{\V^\frac{1}{8}})$ outside of a set with $O(\frac{1}{\V^\frac{1}{4}})$ probability.  As a consequence,
\begin{gather}
|\Tvsd - \bTvsd| = O(\frac{1}{\epsilon \V^\frac{1}{8}}) \\ \notag
\v(\Tvd) = \epsilon^q, \\ \notag
|f - \vs(\Tvd)| \le O(\epsilon^q) \\ \notag
\delta(\Tvd) > 0.
\end{gather}
We define the function $\phi(t) = \int_{\Tvd}^t ds (1-f) - \bp(s)$.  Notice that $\phi$ depends on $\bp(\Tvd)$ although we mostly suppress this dependence.   By the results of Lemma \ref{L:II} we can integrate $\bp$ to arrive at
\begin{equation}
\phi(t) = (1-f)(t - \Tvd) - \frac{1}{\alpha \epsilon} \log(1 + \epsilon \alpha \bp(\Tvd)(t - \Tvd)) + O(\epsilon^q)
\end{equation}
Since $\bp$ is strictly decreasing in Stage II, $\phi$ has a single critical point which we set as $\frac{\se}{\epsilon}$ (as mentioned in section \ref{S:overview}).    We have $\phi'(\frac{\se}{\epsilon}) = (1-f) - \bp(\frac{\se}{\epsilon}) = 0$.  Since we can solve for $\bp$ on Stage $II$, we can solve for $\se$.  In fact, 
\begin{gather}
\se = \epsilon \Tvd + \frac{\delta(\Tvd)}{\alpha(1-f)\bp(\Tvd)} + O(\epsilon^q), \\ \notag
\phi(\frac{\se}{\epsilon}) = \frac{-1}{\alpha \epsilon}\bigg[-\frac{\delta(\Tvd)}{\bp(\Tvd)} + \log(1 + \frac{\delta(\Tvd)}{1-f})\bigg] + O(\epsilon^q), \\ \notag
\phi''(\frac{\se}{\epsilon}) = \epsilon \alpha (1-f)^2 + O(\epsilon^q), \\ \notag
\phi'''(\frac{\se}{\epsilon}) = O(\epsilon^2).
\end{gather}

We split the interval $[\Tvd, \Tvb]$ into three pieces using the times $t_0 = \frac{\se}{\epsilon} - \frac{1}{\epsilon^\kval}$ and $t_1 = \frac{\se}{\epsilon} + \frac{1}{\epsilon^\kval}$ where $\kval$ can taken any value between $\frac{1}{2}$ and $\frac{2}{3}$.  A consequence of our arguments will be  $\Tvd < t_0 < t_1 < \Tvb$.  The following lemma shows that in the time interval $[\Tvd, t_0]$, $\u$ stays close to $\bu$.

\begin{lemma}  \label{L:stochastic_II_part1}
Outside of a set $\Omega$ such that $P(\Omega) \le \epsilon^{q - 2}$ we have for any $t \in [\Tvd, t_0]$,
\begin{gather}  
|\v(t) - \bv(t)| \le \epsilon \bv(t).  \label{E:II_lemma_part1}  \\
|\vs(t) - \bvs(t)| \le \epsilon^\frac{q+1}{2} \\ 
|\p(t) - \bp(t)| \le \epsilon^\frac{q+1}{2}.  \label{E:II_lemma_part1_f}
\end{gather}
\end{lemma}

\begin{proof}
From Lemmas \ref{L:I} and \ref{L:I_and_III_s} we see that (\ref{E:II_lemma_part1})-(\ref{E:II_lemma_part1_f}) hold for $t = \Tvd$.  In Lemma \ref{L:I_and_III_s} we were able to scale $\u - \bu$ by a constant $D$.  Here, things are not so simple.  Set $r = \frac{q+1}{2}$ and define,
\begin{equation}
\w = \bigg(
\begin{array}{c}
\frac{\v - \bv}{h(t)}\\
\frac{\vs - \bvs}{\epsilon^r} \\
\frac{\p - \bp}{\epsilon^r}
\end{array}
\bigg)
\end{equation}
where $h(t) = \bv(\bTvd) \exp[\phi(t)]$.  Note $h(t) = \bv(t)(1 + O(\epsilon^{q-1}))$.  

	We write $\w = (\w_1,\w_2,\w_3)$ and set $T = \inf\{t : \|\w(t)\|_\infty > \epsilon, t < \Tvb \}$.  Now, notice that the coordinates of $a$ given in (\ref{E:def_a}) are all second order polynomials and so an exact second order Taylor expansion exists for $a$.  Using this Taylor expansion gives at time $t \le T$,
\begin{gather}  \label{E:d_w_II}
d\w_1 = g_1(t) \w_1 + O(\epsilon^{r}) + \frac{\sigma_1(u)}{h(t) \sqrt{V}} dB_1 \\
d\w_2 = g_2(t) \w_2 +  O(\epsilon^r) + \frac{\sigma_2(u)}{\epsilon^q \sqrt{V}} dB_2 \\ 
d\w_3 = g_3(t) \w_3 + O(\epsilon^r).
\end{gather}
where
\begin{equation}
g(t) = \bigg(
\begin{array}{c}
-\phi'(t) + (1 - \bvs - \bv - \bp) \\
(f - \bv - 2\bvs) \\
-2\bp\epsilon \\
\end{array}
\bigg)
\end{equation}
We would now like to use $g(t)$ as an integrating factor.  Notice first,
\begin{equation}
g(t) = \bigg(
\begin{array}{c}
O(\epsilon^q) \\
-f + O(\epsilon^q) \\
O(\epsilon) \\
\end{array}
\bigg)
\end{equation}
Set $G_i (s,t) = \exp[\int_s^{t} ds' g_i(s')]$.  Then we can integrate (\ref{E:d_w_II}) using $G$ and arrive at,
\begin{gather}
\w_1(t \wedge T) = \int_{\Tvd}^{t \wedge T} ds G_1(s,t \wedge T)  O(\epsilon^r) 
					+ \int_{\Tvd}^{t \wedge T} dB_1(s) G_1(s,t \wedge T) \frac{\sigma_1(u)}{h(s) \sqrt{V}}  \\
\w_2(t \wedge T) = \int_{\Tvd}^{t \wedge T} ds G_2(s,t \wedge T) O(\epsilon^r) 
					+ \int_{\Tvd}^{t \wedge T} dB_2(s)  G_2(s,t \wedge T) \frac{\sigma_2(u)}{\epsilon^q \sqrt{V}}  \\
\w_3(t \wedge T) = \int_{\Tvd}^{t \wedge T} ds \exp[\int_s^{t \wedge T} ds' g_3(s')] O(\epsilon^{r}).
\end{gather} 
Initial condition terms involving $\w(\Tvd)$ that should appear in the expressions directly above are of lower order, so we have suppressed them for simplicity.  From the above relations we find directly $\w_3(t \wedge T) = O(\epsilon^{r-1})$.  To bound $\w_1$ we consider second moments.
\begin{equation}
E[\sup_{\Tvd \le t' \le t} \w_1^2(t' \wedge T)] \le O(\epsilon^{2(r-1)}) + 
				E[\bigg(\int_{\Tvd}^{t \wedge T} dB(s) \sqrt{\frac{\v(1 - \v - \vs - \p)}{h^2(t)\V}}\bigg)^2]
\end{equation}
where we have used the fact that $G_1(0,t') = O(1)$.  For $t \le T$, we have $\sigma(\v) = O(\v) = O(h(t))$.  Using this observation and standard Martingale arguments gives,
\begin{equation}
E[\sup_{\Tvd \le t' \le t} \w_1^2(t \wedge T)]  
			\le O(\epsilon^{q-1}) +  O(\frac{1}{\bv(\bTvd)} \frac{\int_{\bTvd}^{t \wedge T} ds \exp[-\phi(s)]}{\V})
\end{equation}
The situation is simpler for $\w_2$.  We use $\sigma_2(\u) \le 1$ and arrive at,
\begin{equation}
E[\sup_{t' \le t} \w_2^2(t' \wedge T)] \le O(\epsilon^{2r}) + O(\frac{1}{\epsilon^{2q} \V}).
\end{equation}
Using our moment bounds on $\w_1, \w_2$ and pointwise bound on $\w_3$ we now bound $\w(t)$ and hence remove our restriction of $t \le T$ outside a set of small probability.  Indeed first consider $\w_1$.  By a Chebyshev bound,
\begin{equation}  \label{E:prob_bound}
P(\sup_{\Tvd \le t' \le t_0 \wedge T} |\w_1(t')| \ge \frac{\epsilon}{2}) = O(\epsilon^{q-2}) +  I
\end{equation}
where
\begin{equation}
I = \int_{\Tvd}^t ds \frac{\exp[-\phi(s)]}{\bv(\Tvd) \epsilon^2 \V}
\end{equation}
Recall that $\bv(\bTvd) = \epsilon^q$.  
We consider $I$ by performing a Taylor series expansion of $\phi$ about $\frac{\se}{\epsilon}$,
\begin{align}  \label{E:main_estimate_II}
I & \le \frac{1}{\V \epsilon^{q+2}} \int_{\Tvd}^{t_0} ds \exp[-\phi(s)]
\\ \notag
 & = \frac{1}{\V \epsilon^{q+3}} \int_{\Tvd}^{\se - \epsilon^{1-\kval}} ds \exp[-\phi(\frac{s}{\epsilon})]
\\ \notag
	& = \frac{1}{\V \epsilon^{q+3}} \int_{\epsilon \Tvd}^{\se - \epsilon^{1-\kval}} ds \exp[-\phi(\frac{\se}{\epsilon})]
				O(\exp[\frac{-\alpha (1-f)^2\epsilon}{2}(\frac{s}{\epsilon} - \frac{\se}{\epsilon})^2])
\\ \notag
	& = \frac{\exp[-\phi(\frac{\se}{\epsilon})]}{\V \epsilon^{q+\frac{5}{2}}} 
				O(\exp[-O(\frac{1}{\epsilon^{2\kval}})]).
\end{align}
Now, $\phi$ is a function of $\p(\Tvd)$ for which we have an explicit asymptotic expression.  Using the expansions of Lemma \ref{L:I} we can express $\p(\Tvd)$ in terms of $\p(\TO)$ (the difference is $O(\epsilon |\log \epsilon|)$).  Applying a Taylor series argument for $\phi$, we can then compute $\phi(\frac{\se}{\epsilon})$ in terms of $\p(\TO)$.   Recalling that $\p(\TO) = \frac{1}{1 + \alpha} + O(\epsilon^q)$ we can arrive at the relation between $\philim$ and $\phi$ given in (\ref{E:phi_philim_relation}).  Plugging this into (\ref{E:main_estimate_II}) and using the scaling of Theorem \ref{T:main} gives $I = 	O(\exp[-O(\frac{1}{\epsilon^{2\kval}})])$.  Bounds  for $\w_2, \w_3$ are straightforward and we can conclude $P(T < t_0) \le O(\epsilon^{q-2})$.  
 \end{proof}

Now we consider the interval $[t_0, t_1]$.  This interval is where stochastic effects become important and $\v$ diverges from $\bv$.

\begin{lemma}  \label{L:stochastic_II_part2}
If $\frac{1}{2} < \kval < \frac{2}{3}$ then outside of a set $\Omega$ with $P(\Omega) \le O(\epsilon^2)$ we have,
\begin{gather}
\sup_{t \in [t_0,t_1]} |\p(t) - \bp(t)| < O(\epsilon^\frac{q+1}{2}) \\ \notag
\sup_{t \in [t_0,t_1]}  |\vs(t) - \bvs(t)| > O(\epsilon^\frac{q+1}{2})
\end{gather}
Set,
\begin{equation}
\z_{II} = \sqrt{\frac{1}{\alpha (1-f)^2}} (\frac{\exp[-\phi(\frac{\se}{\epsilon})]}{\V \v(\Tvd) \sqrt{\epsilon}})
\end{equation}
Then,
\begin{equation}  \label{E:z_II_limit}
\lim_{\V \to \infty} \z_{II}  = \sqrt{\frac{1}{\alpha (1-f)^2}} (\frac{\alpha}{1 + \alpha}) \kappa,
\end{equation}
\begin{equation}  \label{E:v_at_t_1}
\v(t_1) = w\big[\sqrt{2\pi}(k+1)\z_{II}\big] \bv(t_0)(1 + O(\epsilon^{2 - 3\kval}),
\end{equation}
and,
\begin{equation}  \label{E:integrated_v}
\int_{t_0}^{t_1} ds \frac{1}{\V \v(s)} = \z_{II} \int_{0}^{\sqrt{2\pi}} ds \frac{1}{w[(k+1)\z_{II} s]}(1 + O(\epsilon^{2 - 3\kval})
\end{equation}
\end{lemma}

\begin{proof}
We first define a stopping time $T$ as follows,
\begin{equation}
T = t_1 \wedge \inf\{t > t_0 : \v(t) > \epsilon^q \text{ or } |\vs(t) - f| \ge \epsilon^q\}.  
\end{equation}
Let $z(t) = \frac{\p - \bp}{\epsilon^r}$ and recall from Lemma \ref{L:stochastic_II_part1} that $r = \frac{q+1}{2}$.  Then $\dot{z}(t \wedge T) = -\epsilon (\p + \bp) z + O(\epsilon^r)$ and we have $z(0) \le O(\epsilon^r)$.  We can conclude $|\p - \bp| = O(\epsilon^r)$ for $t \in [t_0, T]$.  A similar argument shows $|\vs - \bvs| = O(\epsilon^r)$ on $[t_0, T]$.
  
Now we turn to the dynamics of $v(t \wedge T)$.
\begin{equation}  \label{E:first_v_eq}
d\v = \v(\phi'(t) + O(\epsilon^r))dt + \sqrt{\frac{v(k + f + \p(t) + O(\epsilon^q))}{\V}} dB_1
\end{equation}
Next, we linearize $\phi'(t)$ and $\bp(t)$ about $\frac{\se}{\epsilon}$.  Recall $\phi'(\frac{\se}{\epsilon}) = 0$ and $\bp(\frac{\se}{\epsilon}) = (1-f)$.  
\begin{equation}  
d\v = \v(\phi''(\frac{\se}{\epsilon})(t - \frac{\se}{\epsilon}) + O(\epsilon^{2(1-\kval)}))dt + \sqrt{\frac{v(k + 1) + O(\epsilon^{1-\kval}))}{\V}} dB_1
\end{equation}
Now we control $\v$ with the same techniques used in Lemma \ref{L:I_0}.  Define diffusions $\v_A, \v_B$ on the Brownian motion space produced by $B_1,B_2$ as follows,
\begin{equation}
d\v_A = \v\bigg(\phi''(\frac{\se}{\epsilon})(t - \frac{\se}{\epsilon}) + C \epsilon^{2(1-\kval)}\bigg)dt + \sqrt{\frac{\v(k + 1) + O(\epsilon^{1-\kval}))}{\V}} dB_1
\end{equation}
and
\begin{equation}
d\v_B = \v\bigg(\phi''(\frac{\se}{\epsilon})(t - \frac{\se}{\epsilon}) - C \epsilon^{2(1-\kval)}\bigg)dt + \sqrt{\frac{\v(k + 1) + O(\epsilon^{1-\kval}))}{\V}} dB_1
\end{equation}
where $C$ is an $O(1)$ constant.  Then up to time $T$, $\v_B < \v < \v_A$.  We will show that $\v_B(T) \to \v_A(T)$ in distribution and this characterizes $\v(T)$.  Consider $\v_A$.  We proceed as in Lemma \ref{L:I_0}.  First we define,
\begin{equation}
x_A(t) = \frac{\v_A(t)}{\bv(t_0)} \exp[-\phi''(\frac{\se}{\epsilon}) \int_{t_0}^t ds (s - \frac{\se}{\epsilon}) + C \epsilon^{2(1-\kval)}].
\end{equation}
Then define
\begin{equation}
w_A(\tau_A(t)) = x_A(t)
\end{equation}
with
\begin{equation}  \label{E:tau_prime}
\tau_A'(t) =   \frac{k+1+O(\epsilon^{1-\kval})}{\V \v(t_0)}
	\exp[-\phi''(\frac{\se}{\epsilon}) \int_{t_0}^t ds (s - \frac{\se}{\epsilon}) + C\epsilon^{2(1-\kval)}]
\end{equation}

Under these transformations, $w_A$ is a weak solution of the following Feller diffusion.
\begin{equation}  \label{E:Feller_w}
dw_A = \sqrt{w_A} dB_1.
\end{equation}
We would like to integrate (\ref{E:tau_prime}) to obtain $\tau_A$.  This will be made simpler with the following formula for $\bv(t_0)$ in terms of $\bv(\Tvd)$.   
\begin{align}  \label{E:t_0_to_T_I}
\bv(t_0) & = \bv(\Tvd) \exp[\phi(t_0)] (1 + O(\epsilon))
\\ \notag
	& = \bv(\Tvd) \exp[\phi(\frac{\se}{\epsilon}) - \phi''(\frac{\se}{\epsilon})\int_{t_0}^\frac{\se}{\epsilon} ds (s - \frac{\se}{\epsilon}) + O(\epsilon^{2 - 3\kval}))].
\end{align}
Integrating (\ref{E:tau_prime}) and applying the above relation gives the following,
\begin{align}  \label{E:tau_A_t_1}
\tau_A(t) 
& = \exp[-\phi(\frac{\se}{\epsilon})] \int_{t_0}^{t} dt \frac{k+1+O(\epsilon^{1-\kval})}{\V \bv(\Tvd)}
\exp[\phi''(\frac{\se}{\epsilon}) \int_{\frac{\se}{\epsilon}}^{t}ds (s - \frac{\se}{\epsilon}) + O(\epsilon^{2 - 3\kval})]
\\ \notag
& = (\frac{\exp[-\phi(\frac{\se}{\epsilon})](k+1))}{\bv(\Tvd)\V \sqrt{\epsilon}})(\sqrt{\frac{1}{\alpha(1-f)^2}})
\Phi\bigg(\sqrt{\epsilon \alpha (1-f)^2}(t - \se)\bigg)(1 + O(\epsilon^{2 - 3\kval}))
\end{align}
where $\Phi(t) = \int_{-\infty}^t ds \exp[-\frac{s^2}{2}]$.

	With (\ref{E:t_0_to_T_I}) and (\ref{E:tau_A_t_1}) we can compute $\v_A(t)$.  We have,
\begin{align}  \label{E:final_vA_formula}
\v_A(t) 
	& = w\big[\tau_A(t)\big] \v(t_0) \exp[\phi''(\se) \int_{t_0}^{t} ds (s - \frac{\se}{\epsilon}) + O(\epsilon^{2-3\kval})]
\\ \notag
	& = w\big[\tau_A(t)\big] \bv(t_0) \exp[\epsilon \alpha (1-f)^2 \int_{t_0}^t ds (s - \frac{\se}{\epsilon})](1 + O(\epsilon^{2 - 3\kval})).
\end{align}
The results for $v_B$ are identical and so we can conclude $\v(t \wedge T) = \v_A(t \wedge T)(1 + O(\epsilon^{2-3\kval)})$.  As described in the proof of Lemma \ref{L:stochastic_II_part1}, a Taylor series argument on $\phi$ using the expansions for $\bp$ developed in Lemma \ref{L:I} allows us to derive (\ref{E:phi_philim_relation}).   Using this formula and our scaling for $\epsilon$ gives (\ref{E:z_II_limit}).   

 Now we eliminate the $T$ dependence of our result.  Since $\v_A(t) = \bv(t_0) w[\tau(t)]$ and $\tau(t) = O(1)$ we can conclude that $\v_A(t) \ll \epsilon^q$ with high probability.  More precisely, through a Chebyshev inequality we have
\begin{equation}
P(\sup_{t \in [t_0,t_1]} \v_A(t) > \epsilon^q) \le P(w\big[\tau(t)\big] > \frac{\epsilon^q}{\bv(t_0)}) \le P(w[\tau(t)] > O(\exp[-O(\frac{1}{\epsilon^\kval})]))
	\le O(\epsilon^2).
\end{equation}

(\ref{E:v_at_t_1}) now follows by plugging in $t=t_1$ in (\ref{E:final_vA_formula}), and (\ref{E:integrated_v}) follows by using  (\ref{E:final_vA_formula}) in the integral $\int_{t_0}^{t_1} ds \frac{1}{\V \v(s)}$ and applying the substitution, $s \to \tau_A(s)$.  
 \end{proof}

	Finally assuming that $\v(t_1) \ne 0$, we consider the dynamics of $\v$ on $[t_1, \Tvb]$.  
	
\begin{lemma}  \label{L:stochastic_II_part3}
Assume $\v(t_1) = \eta \v(t_0)$ for some $\eta > 0$.    Then for $t \in [t_1, \Tvb]$ outside a set $\Omega$ with $P(\Omega) < \epsilon^{q-2}$ we have,
\begin{gather}
|\vs(t) - f| \le O(\epsilon^\frac{q+1}{2})  \\ 
\p(\Tvb) = \bp(\bTvb) - \epsilon \alpha^2 \frac{\bp^2(\bTvd)}{1 - f - \bp(\bTvb)} \log \eta + O(\epsilon^2).  \label{E:p_stoch_perturb}
\end{gather}
\end{lemma}

\begin{proof}
After time $t_1$ the system returns to deterministic behavior. By the same arguments used in Lemma \ref{L:stochastic_II_part1} we have,
\begin{gather}
|\vs(t) - \bvs_1(t)| \le \epsilon^\frac{q+1}{2} \\ 
|\p(t) - \bp_1(t)| \le \epsilon^\frac{q+1}{2}.
\end{gather}

In Lemma \ref{L:stochastic_II_part1} we showed that in $[\Tvd, t_0]$, $\v$ is well approximated by $\bv$.  The same holds in $[t_1, \Tvb]$, except that now we must restart the deterministic system so that $\bv(t_1) = \v(t_1)$.  With this in mind we can apply the arguments of Lemma \ref{L:stochastic_II_part1} to justify the following relation.
\begin{equation}
\v(t)  = \eta \v(t_1) \exp[\phi(t) - \phi(t_1) + O(\epsilon^\frac{q+1}{2})] 
\end{equation}
Then if we use the arguments of Lemma \ref{L:stochastic_II_part1} in which we prove $|\v - \bv| \le \epsilon \bv$ we can argue as follows,
\begin{align}
\v(t) 
	& = \eta \bv(t) (1 + O(\epsilon)) \exp[-(\phi(t_1) - \phi(t_0))]
\\ \notag
	& = \eta \bv(t)(1 + O(\epsilon)) \exp[-\bigg(\phi''(\se)((t_1 - \se)^2 - (\se - t_0)^2)\bigg) + O(\bar{\phi}'''(\se)(t_1 - t_0)^3)]
\\ \notag
	& = \eta \bv(t)(1 + O(\epsilon^{2 - 3\kval}))
\end{align}
Now we consider $\Tvb$ in comparison to $\bTvb$.  Perturbing off of $\bTvb$ gives,
\begin{equation}
\v(\bTvb + \Delta t)
	= \eta \bv(\bTvb)(1 + O(\epsilon ^{2 - 3\kval}))
	\exp[(1-f)\Delta t - \frac{1}{\alpha \epsilon} \log(1 + \alpha \epsilon \bp(\bTvb)) \Delta t)]
\end{equation}
A Taylor expansion on $\Delta t$ leads to
\begin{align}  \label{E:taylor_delta_t}
\v(\bTvb + \Delta t)
	& = \eta \bv(\bTvb)(1 + O(\epsilon ^{2 - 3k}))
	\exp[(1-f)\Delta t  - \bp(\bTvb) \Delta t + O(\Delta t^2)\bigg)]
\end{align}
We want to find $\Delta t$ such that $\v(\bTvb + \Delta t) = \epsilon^q$ since then we will have $\Tvb = \bTvb + \Delta t$.  Solving using (\ref{E:taylor_delta_t}) gives,
\begin{equation}
\Delta t = -(\frac{1}{\delta(\bTvb)}) \log \eta + O(\epsilon^2)
\end{equation}
Our main interest in determining $\Tvb$ is our need to compute an expansion for $\p(\Tvb)$.  We have,
\begin{align}
\p(\Tvb)
	& = \bp(\bTvb) + \bp'(\bTvb) \frac{\Delta t}{\epsilon} + O(\epsilon^2)
\\ \notag
	& = \bp(\bTvb) - \epsilon 
	\frac{\bp^2(\bTvb)}{\delta(\bTvb)}) \log \eta + O(\epsilon^2)
\end{align}
 \end{proof}

\subsection{Stage IV}

Stage IV is similar to Stage II.  As in Stage II, using Lemmas \ref{L:III} and \ref{L:I_and_III_s} we have outside of a set of vanishing probability,
\begin{gather}
|\Tvsd - \bTvsd| \le O(\frac{1}{\epsilon \V^\frac{1}{4}}) \\ \notag
|\v(\Tvsd) - (1 - p(\Tvsd))| = O(\epsilon), \\ \notag
\vs(\Tvsd) = \epsilon^q \\ \notag
\delta(\Tvsd) < 0.
\end{gather}
We set $\psi(t) = \int_{0}^t ds \bp(s) - (1-f)$.  Then $\psi$ has a single critical point which we label $\frac{\ses}{\epsilon}$.  In Stage II we define $t_0 = \frac{\se}{\epsilon} - \frac{1}{\epsilon^\kval}$ and $t_1 = \frac{\se}{\epsilon} + \frac{1}{\epsilon^\kval}$.  In this section, for Stage IV we define $t_0^* = \frac{\ses}{\epsilon} - \frac{1}{\epsilon^\kval}$ and $t_1^* = \frac{\ses}{\epsilon} + \frac{1}{\epsilon^\kval}$. 

	The proofs of Stage IV are almost identical to those of Stage II.  We state the analogues of Lemmas \ref{L:stochastic_II_part1} and \ref{L:stochastic_II_part3} since the proofs follow identical arguments.  In Lemma \ref{L:stochastic_IV_part3} we only keep $O(1)$ terms for the $\p(\Tvsb)$ because we do not need the $O(\epsilon)$ accuracy in the later cycles that we need in the first cycles.

\begin{lemma}  \label{L:stochastic_IV_part1}
Outside of a set $\Omega$ such that $P(\Omega) \le O(\epsilon^2)$ we have for any $t \in [\Tvsd, t_0^*]$,
\begin{gather}
|\v(t) - \bv(t)| \le \epsilon^\frac{q+1}{2}. \\
|\vs(t) - \bvs(t)| \le \epsilon \bvs(t) \\ 
|\p(t) - \bp(t)| \le \epsilon^\frac{q+1}{2}.
\end{gather}
\end{lemma}

\begin{lemma}  \label{L:stochastic_IV_part3}
Suppose $\vs(t_1^*) = \eta \vs(t_0^*)$.  Then for $t \in [t_1, \Tvsb]$,  
\begin{gather}
|\v(t) - (1-p(t))| \le O(\epsilon)  \\ \notag
\p(\Tvsb) = \bp(\bTvsb) + O(\epsilon)
\end{gather}
\end{lemma}

	Now we turn to the analogue of Lemma \ref{L:stochastic_II_part2}.  The arguments are essentially the same, but we consider $\psi$ instead of $\phi$.  We have,
\begin{gather}  \label{E:psi_ses_values}
\ses  = \epsilon \Tvsd  + \log[\frac{(1-f)(1 - (1 + \alpha)\p(\Tvsd))}{(f - \alpha(1-f))\p(\Tvsd)}], \\ \notag
\psi(\frac{\ses}{\epsilon}) = \frac{1}{\epsilon}
	\bigg[\frac{1}{1+\alpha} \log(\frac{1 - (1+\alpha)\p(\Tvsd)}{f - \alpha(1-f)})
	 - (1-f)\log(\frac{(1-f)(1 - (1+\alpha)\p(\Tvsd))}{(f - \alpha(1-f))\p(\Tvsd)})\bigg] \\ \notag
\psi''(\frac{\ses}{\epsilon}) = \epsilon(f - \alpha(1-f)) + O(\epsilon^2). \\ \notag
\psi'''(\frac{\ses}{\epsilon}) = O(\epsilon^2).
\end{gather}
Then the analogue to Lemma \ref{L:stochastic_II_part2} is the following result.  We do not need the expression $\int_{t_0^*}^{t_1^*} ds \frac{1}{\V \vs(s)}$ because the lineages of the mutant type must coalesce to the original mutant cell by time zero.

\begin{lemma}  \label{L:stochastic_IV_part2}
If $\frac{1}{2} < \kval < \frac{2}{3}$ then outside of a set $\Omega$ with $P(\Omega) \le O(\epsilon^2)$ we have,
\begin{gather}
\sup_{t \in [t_0,t_1]} |\p(t) - \bp(t)| < O(\epsilon^\frac{q+1}{2}) \\ \notag
\sup_{t \in [t_0,t_1]}  |\vs(t) - \bvs(t)| > O(\epsilon^\frac{q+1}{2})
\end{gather}
Set,
\begin{equation}
\z_{IV} = \sqrt{\frac{1}{(1-f)(f - \alpha(1-f))}} (\frac{\exp[-\psi(\frac{\ses}{\epsilon})]}{\V \bv(\bTvsd) \sqrt{\epsilon}})
\end{equation}
Then,
\begin{equation} \label{E:z_IV_limit}
\lim_{\V \to \infty} \z_{IV} = \sqrt{\frac{1}{(1-f)(f-\alpha(1-f))}}(\frac{1}{f})\bigg(\frac{1 + \alpha(1+H)}{(1+\alpha)\alpha(1+H)}\bigg)^H (\eta_{IV})^\frac{\alpha}{H} \kappa.
\end{equation}
(recall the definition of $\eta_{IV}$ in (\ref{E:eta_IV_def}) from Theorem \ref{T:main}),
\begin{equation} 
\vs(t_1) = w\big[\sqrt{2\pi}(k^*+f)\z_{IV}\big] (\frac{1 - \bp(\bTvsd)}{f}) \exp[-\delta(\bTvsd)] \bvs(t_0)(1 + O(\epsilon^{2 - 3\kval}))
\end{equation}
\end{lemma}

	The expression $\eta_{IV}$ in (\ref{E:z_IV_limit}) requires explanation.  If the wild type is not lost, then it experiences a stochastic perturbation in Stage II.  As Lemmas \ref{L:stochastic_II_part2} and \ref{L:stochastic_II_part3} show, this perturbation is $O(1)$ and influences $\p(\Tvb)$ by $O(\epsilon)$.  The $O(\epsilon)$ perturbation on $\p(\Tvb)$ is integrated over Stage IV which is of duration $O(\frac{1}{\epsilon})$ and so the perturbation has an $O(1)$ effect on $\vs$ during Stage IV.  This is where $\eta_{IV}$ comes from.  More specifically, when we expand $\psi$ in Taylor series to obtain an expansion in terms of $\psilim$, the $O(\epsilon)$ term in (\ref{E:p_stoch_perturb}) in Lemma \ref{L:stochastic_II_part3} is responsible for $\eta_{IV}$.

\appendix
\section*{Appendix}
\renewcommand{\theequation}{A.\arabic{equation}}
\renewcommand{\thefigure}{A.\arabic{figure}}
\setcounter{equation}{0}

	Here, we explain the rescaling of section \ref{S:bd} precisely.  We introduce dimensional constants $\V, \P, \T$ for the units of infected cells, CD8 cells, and time respectively and define the non-dimensional variables $\tilde{\v}, \tilde{\vs}, \tilde{\p}$ by setting
\begin{equation}
\tilde{\v}(\tilde{t}) = \frac{\v(\T\tilde{t})}{\V}, \quad
\tilde{\vs}(\tilde{t}) = \frac{\vs(\T\tilde{t})}{\V}, \quad
\tilde{\p}(\tilde{t}) = \frac{\p(\T\tilde{t})}{\P}, \quad
\tilde{t} = \frac{t}{\T}.
\end{equation}
We can then hope to approximate the birth-death process through the following SDE.
\begin{gather}  \label{E:unscaled_SDE}
d\tilde{\v} = \T \tilde{\v}(\Delta k - c\V(\tilde{\v} + \tilde{\vs}) - a\P\tilde{\p}) d\tilde{t}
	+ \sqrt{\frac{\tilde{\v}\T(k + c\V(\tilde{\v} + \tilde{\vs}) + a\P\tilde{\p})}{\V}} dB_1(\tilde{t}),
\\ \notag
d\tilde{\vs} = \T \tilde{\vs}(\Delta k^* - c\V(\tilde{\v} + \tilde{\vs})) d\tilde{t}
	+ \sqrt{\frac{\tilde{\vs}\T(k^* + c\V(\tilde{\v} + \tilde{\vs}))}{\V}} dB_2(\tilde{t}), 
\\ \notag
d\tilde{\p} = \T \tilde{p}(b\V\tilde{v} - d\P\tilde{p}) d\tilde{t}
	+ \sqrt{\frac{\tilde{p}\T(h + b\V\tilde{\v} + d\P\tilde{p})}{\P}} dB_3(\tilde{t}).
\end{gather}
We choose $\T$, $\V$, and $\P$ so that $\Delta k \T = 1$, $c \V \T = 1$, and $a \P \T = 1$.  Plugging this into (\ref{E:unscaled_SDE}) gives
\begin{gather}  \label{E:rescaled_SDE}
d\tilde{\v} =  \tilde{\v}(1 - (\tilde{\v} + \tilde{\vs}) - \tilde{\p}) d\tilde{t}
	+ \sqrt{\frac{\tilde{\v}(\tilde{k} + (\tilde{\v} + \tilde{\vs}) + \tilde{\p})}{\V}} dB_1(\tilde{t}),
\\ \notag
d\tilde{\vs} = \tilde{\vs}(f - (\tilde{\v} + \tilde{\vs})) d\tilde{t}
	+ \sqrt{\frac{\tilde{\vs}(\tilde{k^*} + (\tilde{\v} + \tilde{\vs}))}{\V}} dB_2(\tilde{t}), 
\\ \notag
d\tilde{\p} = \epsilon \tilde{p}(\tilde{v} - \alpha \tilde{p}) d\tilde{t}
	+ \sqrt{\frac{\epsilon \tilde{p}(\tilde{h} + \tilde{\v} + \alpha\tilde{p})}{\P}} dB_3(\tilde{t}).
\end{gather}
where,
\begin{equation}
\tilde{k} = \frac{k}{\Delta k}, \quad \tilde{k^*} = \frac{k^*}{\Delta k}, \quad \epsilon = \frac{b}{c},  \quad
\alpha = \frac{d}{a}.
\end{equation}
We assume that the coefficients in (\ref{E:rescaled_SDE}) are all $O(1)$.  This will be true if $\V, \P$ are on the order of the infected cell and CTL population counts and the system is assumed to vary on that scale.  We now drop the tildes and, for simplicity, the variance terms in the $\tilde{\p}$ equation.  This gives (\ref{E:SDE})


\newcommand{\noopsort}[1]{} \newcommand{\printfirst}[2]{#1}
  \newcommand{\singleletter}[1]{#1} \newcommand{\switchargs}[2]{#2#1}

\end{document}